\documentclass{article}
\usepackage[utf8]{inputenc}
\usepackage{amsmath,amsfonts,amssymb,graphicx,algorithm,caption}
\usepackage{cite}
\usepackage{mathtools}
\usepackage{booktabs}
\usepackage[margin=.7in]{geometry}
\usepackage{tabularx}
\usepackage{graphicx} 
\usepackage[utf8]{inputenc}
\usepackage{enumitem}
\usepackage[english]{babel}
\usepackage{url}
\usepackage{hyperref}
\usepackage{cleveref}
\usepackage{amsthm}

\newtheorem{theorem}{Theorem}
\newtheorem{lemma}[theorem]{Lemma}
\newtheorem{corollary}[theorem]{Corollary}
\newtheorem{conjecture}[theorem]{Conjecture}

\theoremstyle{definition}

\newtheorem{observation}[theorem]{Observation}
\newtheorem{rem}[theorem]{Remark}
\newtheorem{example}[theorem]{Example}

\newcommand{\mex}{\mathrm{mex}}
\renewcommand{\P}{\ensuremath{\mathcal P}}
\setlist[itemize]{ topsep=0pt}
\begin{document}

\title{Additive Sink Subtraction}

\author{
Anjali Bhagat\thanks{Department of Industrial Engineering and Operations Research, IIT Bombay, India}
\and
Urban Larsson\thanks{Department of Industrial Engineering and Operations Research, IIT Bombay, India. Corresponding author: larsson@iitb.ac.in}
\and
Hikaru Manabe\thanks{College of Information Science, School of Informatics, University of Tsukuba, Japan}
\and
Takahiro Yamashita\thanks{Department of Mathematics, Hiroshima University, Japan. Funding: JST SPRING, Grant Number JPMJSP2132}
}

%\date{}

\maketitle

\begin{abstract}
Subtraction games are a classical topic in Combinatorial Game Theory. An early result of Golomb~(1966) shows that every subtraction game with a finite move set has an eventually periodic nim-sequence, but the known proof yields only an exponential upper bound on the period length. Even three-move games are not yet fully understood. Nevertheless, Flammenkamp (1997) conjectures a striking classification: non-additive rulesets have  linear period lengths of the form ``the sum of two moves'', where the choice of which two moves displays fractal-like behavior, while additive sets $S=\{a,b,a+b\}$ have purely periodic outcomes with linear or quadratic period lengths. Despite early attention in \emph{Winning Ways} (1982), the general additive case remains open. We analyze a dual winning convention, called {\sc sink subtraction}. Unlike the standard wall convention, where moves to negative positions are forbidden, the sink convention declares a player the winner upon moving to a non-positive position. We show that {\sc additive sink subtraction} admits a complete solution: its nim-sequence is purely periodic with an explicit linear or quadratic period formula. We further conjecture a nimber duality between additive sink subtraction and classical wall subtraction. 

\end{abstract}

\section{Introduction}
\begin{figure}
    \centering
    \includegraphics[
        width=0.3\linewidth,
        trim=0cm 8cm 0cm 5.5cm,
        clip
    ]{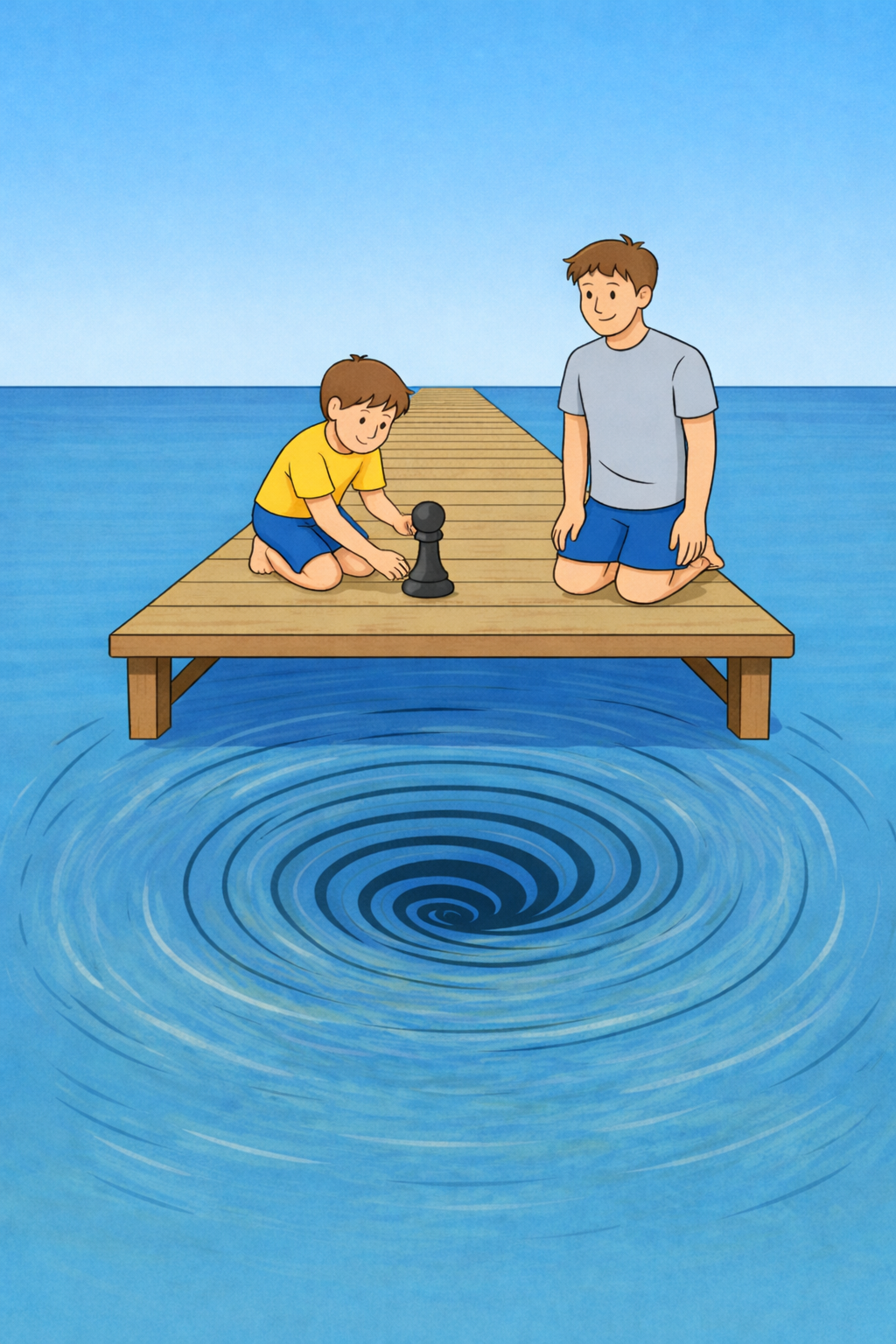}
    \includegraphics[
        width=0.32\linewidth,
        trim=0cm 11cm 0cm 5cm,
        clip
    ]{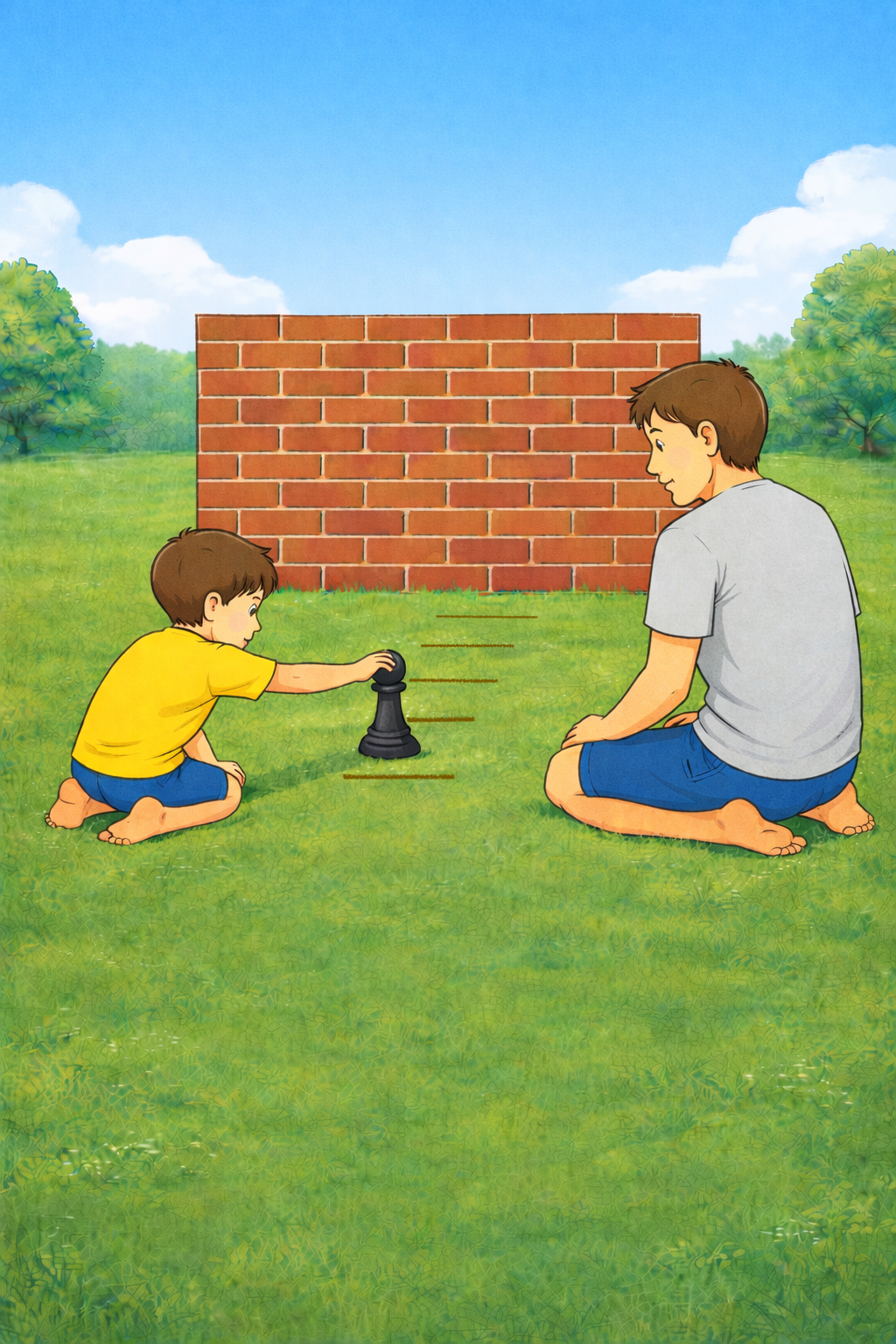}
    \caption{{\sc sink} and {\sc wall subtraction}, respectively. Suppose the ruleset is $S=\{1,5\}$. To the Left, the boy can choose to move into the sink or he can play to the tip of the jetty, and let daddy win. To the Right, the boy has only one move, because the larger move is prohibited by the wall. Here, the wall game is determined already by parity.}
    \label{fig:sinkwall}
\end{figure}

We study classical combinatorial subtraction games \cite{berlekamp2004winning}, with a twist on the terminal condition. Let $S\subset \mathbb N=\{1,2,\ldots \}$. For any $x\in \mathbb N$ and all $s\in S$, $x-s\in\mathbb Z $ is a move option. 
The {\em sink} is an infinite sequence of terminal $\mathcal P$-positions ($\mathcal{P}$revious player winning positions) for all non-positive integers. No move is possible from a position in the sink. We study two-player impartial games in the normal-play convention, so ``a player who cannot move loses'', and this is the same as ``a player who plays into the sink wins''. We denote this ruleset family by {\sc sink subtraction}. Observe that, if $S =\mathbb N$, we play the game of {\sc nim} and {\sc sink subtraction} coincides with {\sc wall subtraction} (our terminology), where all negative heap sizes are forbidden. In Figure~\ref{fig:sinkwall}, we illustrate with the ruleset $S=\{1,5\}$ how sink vs. wall play can differ.

\begin{rem}
In the classical normal-play setting, we may imagine that a concrete wall is placed at position $-1$, preventing players from entering any negative number. Curiously, if we force negative positions to be terminal N-positions, information is lacking about how to play a disjunctive sum of games; we need to assign non-zero nimbers to every nonpositive position. Hence, in this sense, the standard wall convention actually represents infinitely many possible termination rules, while the sink convention yields a unique interpretation (since all P-positions have nim-value zero). 
\end{rem}

Indeed, if $S$ is finite, the two ruleset families {\sc sink} and {\sc wall} usually satisfy distinct {\em nimber} patterns. A game's nimber (a.k.a. nim-value) is defined via the standard minimum exclusive rule; for any finite $X\subset \mathbb N_0$, let $\mex X = \min \{\mathbb N_0\setminus X\}$, where $\mathbb N_0 = \mathbb N\cup \{0\}$. Then the nim-value algorithm of $x$ is $v(x)=\mex \{v(x')\}$, where $x'$ ranges over the options of $x$. Note that this gives value zero for all positions in the sink.

For example, if $S=\{2,5\}$, then starting from $x=0$, {\sc wall} has a purely periodic nim-sequence of period length 7, $${\bf 0,0,1,1,0,2,1},0,0,1,1,0,2,1\ldots, $$ while {\sc sink}, starting from $x=1$, begins $$1,1,2,{\bf 2,1, 0,0,1,1,0},2,1,0,0,1,1,\ldots$$ with a pre-period of size three but the same period length 7. The periodic patterns are the same but occur in a different order: 0011021 and 2100110 respectively. 

To the best of our knowledge, {\sc sink subtraction} has only been studied in the literature in terms of outcomes, and only in one paper \cite{althofer1995superlinear}. They remark that one can relate the outcomes of the two conventions via the formula 
\begin{align}\label{eq:outwallsink}
o_{\rm wall}(x)=o_{\rm sink}(x+\max(S)+1). 
\end{align}
The proof of this boils down to the observation that the first $\max S$ positions in the sink convention are $\mathcal N$-positions (curre$\mathcal{N}$t player winning positions). Hence, onwards, the shifted outcome patterns must be the same. However, when we study nimbers, the implication is not this direct, and one has to be more careful. In our above example, we can replace outcome with nim-value, and the formula still holds. However, as we will see in Example~\ref{ex:257}, this is not always true. %In case of {\sc additive subtraction}, an alternative formula is conjectured  just before that example.  

In this paper we continue the pursuit of \cite{althofer1995superlinear}, by studying the more general nim-values of {\sc sink subtraction}. We use the convention of position $x=0$ belonging to the sink, so we start nim-sequences by the value (one) of position one. 
 Of course, the folklore periodicity theorem holds for rulesets where the terminal condition is a sink instead of a `wall', with exactly the same proof. 
\begin{theorem}[Folklore \cite{golomb1966mathematical}]\label{thm:folk}
    Any instance of {\sc sink subtraction} has an ultimately periodic nim-sequence.
\end{theorem}
\begin{proof}
Let $S$ be the subtraction set of {\sc sink subtraction}.  %Let $\max S \in S$ be the largest element in the set $S$.
The {\em window} of options for any position in a subtraction game is determined by $\max S$. The window symbols are nimbers, computed recursively. By the mex-rule, the number of nimbers is bounded by $\rho = |S|$.

%Through the $\mathrm{mex}$ (minimum excluded) rule to determine the nimber value of any impartial game, we know that any position in these games cannot exceed $\max S$. 

Thus, the number of possible combinations of symbols in a window of size $\max S$ is $\rho ^{\max S}$. This number is finite as $\max S$ is finite. 

The window of options of a position $x$ is $\{x-\max S,\ldots ,x-1\}$, while for the position $x+1$, the window of options is $\{x+1-\max S,\ldots , x\}$, and so on, thus making the number of windows of options infinite. 

The number of combinations of symbols is finite, while the number of windows of options is infinite; this implies that some window symbol combination will repeat. The first occurrence of such a repetition will define the period length. It is possible that, for an initial sequence, window symbol combinations do not repeat; if so, this forms the pre-period. %After this we will witness expansion, and it will be periodic. Also, it may happen that there is no pre-period for some games.
\end{proof}

Note that Theorem~\ref{thm:folk} establishes eventual periodicity for any terminating condition; in general, one can assign arbitrary patterns of forced nimbers for negative integers, and the folklore theorem still holds. There is a recent paper that discusses such settings \cite{miklós2023superpolynomial}, and they establish super-polynomial period length for three-move games with specific such relaxed winning conditions. 

%The {\em total expansion period length} is the period length of the expanded array over all sink sizes.
%\begin{conjecture}
%  For a fixed subtraction set, regardless of cycle length, the total expansion period length remains the same.
%\end{conjecture}

The purpose of this paper is to establish the period lengths of the ruleset family {\sc additive sink subtraction} $S$ with subtraction set 
\begin{align}\label{eq:S}
S=\{s_1,s_2,s_1+s_2\},
\end{align}
where $s_1,s_2\in \mathbb N$, $s_1<s_2$. There is a conjecture from Flammenkamp's PhD thesis \cite{Flammenkamp_1997} concerning the periodicity of  {\sc additive wall subtraction}, and this partially paraphrases the treatment in Winning Ways \cite{berlekamp2004winning}, where the linear-period-length part of the conjecture has been explained (without a proof). As far as the authors know, the non-linear (quadratic) part has not yet been fully solved in the literature, but conjectured also in \cite{ward2016conjecture}.\footnote{Non-additive three-move games are conjectured to have period lengths of the form ``the sum of two moves'', where ``which two moves?'' appear to satisfy some fractal patterns \cite{Flammenkamp_1997}.} Recently a subset of the authors completely determine the closed formulas of the nimbers in the {\em primitive quadratic regime} \cite{LM} (whenever $a<\delta <2a$ and $\gcd{a,\delta}=1$), by applying elementary combinatorial number theory to the \P-position bracket expression from \cite{berlekamp2004winning} (see also Section~\ref{sec:final}). 

Let us state our main result. In the statement we use a different representation of \eqref{eq:S}. 

%We take $d$ as $ d\equiv \delta\ \mathrm{mod}(2m)$. We find that the period length is related to $m,2m+\delta$ and \\ $j=m/\mathrm{gcd}(m-k,k), \text{ where }k= \delta\pmod m$. 

\begin{theorem}[Main Theorem]\label{thm:main} 
Let $m,\delta\in \mathbb N$.  
    Consider {\sc additive sink subtraction} with subtraction set $S(m,\delta)= \{m,m+\delta,2m+\delta\}$, and let $ d := \delta \ \pmod {2m} $, with $0 \le d < 2m$. Then $S$ has purely periodic nim-values, with period length  
\[ p(m,\delta) = \begin{cases}
        3m+2\delta-d, & \mathrm{if}\ d\in [0,m];\\ 
        m(m+2\delta+d)/\gcd(m,d), & \mathrm{otherwise }.
        % m(3m+2\delta-(2m-d))/\mathrm{gcd}(m-d_1,d_1); & \mathrm{otherwise}.
       \end{cases}
\]
\end{theorem}
The immediate consequence for the outcomes of {\sc additive wall subtraction} is as follows. 
\begin{corollary}\label{cor:add} 
Let $m,\delta\in \mathbb N$, and consider {\sc additive wall subtraction} with subtraction set $S(m,\delta)= \{m,m+\delta,2m+\delta\}$. Then the period length of the outcomes is the same as in Theorem~\ref{thm:main}, $p(m,\delta)$.
\end{corollary}
\begin{proof}
    This is obvious by~\eqref{eq:outwallsink}.
\end{proof}
We conjecture that there is an explicit duality formula for relating the nimbers of {\sc additive sink subtraction} with those of {\sc additive wall subtraction}, of the form, for all $x$, 
\begin {align}\label{eq:conj}
v_{\rm wall}\bigl ([x]_p\bigr) = \sigma\bigl(v_{\rm sink}\bigl([-s_2-x]_p\bigr)\bigr), 
\end{align}

where $\sigma : 0\rightarrow 2, 1\rightarrow 1, 2\rightarrow 0, 3\rightarrow 3$.  Here $[x]_p$ denotes the heap size modulo the (outcome) period length $p$. 
\iffalse
\begin{conjecture}
Let $S=\{s_1, s_2, s_3\}$ be a subtraction set satisfying $s_1<s_2, s_1+s_2=s_3, s_1, s_2, s_3 \in \mathbb{N}$, and let $\sigma : 0\rightarrow 2, 1\rightarrow 1, 2\rightarrow 0, 3\rightarrow 3$. Let $v_{\rm wall}(x)$ and $v_{\rm sink}(x)$ denote the nim-values of the wall and sink subtraction games with set $S$ at position $x$, respectively. We write $G \to H$ to indicate that position $H$ is reachable from $G$ in a single move. 
Recall, if $x'-x\in S$, then $x$ is an option of $x'$ [[use this style instead?]]. 
Then, the following properties hold:
\begin{enumerate}
\item[(i)] $v_{\rm wall}(x) = \mathrm{mex}({v_{\rm wall}(x'): x' \rightarrow x}); $
\item[(ii)] $\sigma(\mathrm{mex}(T)) = \mathrm{mex}(\sigma(T))$ for $T = \{v_{\rm wall}(x'): x\rightarrow x'\}$;
\item[(iii)] $v_{\rm wall}(x) = \sigma(v_{\rm sink}(2p-s_2-x+1))$,
\end{enumerate}
where $p$ is as in Theorem~\ref{thm:main}. 
\end{conjecture}
\fi

\begin{example}\label{ex:257}
The {\sc wall} and {\sc sink} nim-sequences of $S=\{2, 5, 7\}$ are displayed in the below table; they have period length $22$. The upper (lower) row represents the heap sizes for {\sc wall} ({\sc sink}), and the middle row entries illustrate the respective nimbers ({\sc wall} upper and {\sc sink} lower).  Note that the values satisfy Equation~\ref{eq:conj}.

\iffalse
The pattern for $v_{\rm wall}(x)$ is below for $0\leq x\leq p$.
\begin{table}[htbp!]
\def\arraystretch{1.2}
\begin{tabular}{|c|c|c|c|c|c|c|c|c|c|c|c|c|c|c|c|c|c|c|c|c|c|c|}
\hline
$x$ & $0$ & $1$ & $2$ & $3$ & $4$ & $5$ & $6$ & $7$ & $8$ & $9$ & $10$ & $11$ & $12$ & $13$ & $14$ & $15$ & $16$ & $17$ & $18$ & $19$ & $20$ & $21$ \\ \hline
$v_{\rm wall}(x)$ & $0$ & $0$ & $1$ & $1$ & $0$ & $2$ & $1$ & $3$ & $2$ & $2$ & $0$ & $3$ & $1$ & $0$ & $0$ & $1$ & $1$ & $2$ & $2$ & $3$ & $3$ & $2$  \\ \hline
\end{tabular}
\vspace{3mm}
\caption{TBA}
\label{tab:wall_257_nimbers}
\end{table}
$v_{\rm wall}(x) = 0, 0, 1, 1, 0, 2, 1, 3, 2, 2, 0, 3, 1, 0, 0, 1, 1, 2, 2, 3, 3, 2$
The pattern for $v_{\rm sink}(x)$ is below for $1\leq x\leq p+1$.
\begin{table}[htbp!]
\def\arraystretch{1.2}
\begin{tabular}{|c|c|c|c|c|c|c|c|c|c|c|c|c|c|c|c|c|c|c|c|c|c|c|}
\hline
$x$ & $1$ & $2$ & $3$ & $4$ & $5$ & $6$ & $7$ & $8$ & $9$ & $10$ & $11$ & $12$ & $13$ & $14$ & $15$ & $16$ & $17$ & $18$ & $19$ & $20$ & $21$ & $22$ \\ \hline
$v_{\rm sink}(x)$ & $1$ & $1$ & $2$ & $2$ & $1$ & $3$ & $2$ & $0$ & $0$ & $3$ & $1$ & $0$ & $2$ & $1$ & $1$ & $2$ & $2$ & $0$ & $3$ & $3$ & $0$ & $0$  \\ \hline
\end{tabular}
\vspace{3mm}
\caption{TBA}
\label{tab:sink_257_nimbers}
\end{table}
$v_{\rm sink}(x) = 1, 1, 2, 2, 1, 3, 2, 0, 0, 3, 1, 0, 2, 1, 1, 2, 2, 0, 3, 3, 0, 0$
\fi
\begin{table}[htbp!]
%\caption{The middle row illustrates the nimbers at the heap sizes for the wall and sink convention respectively.}
%\label{tab:wall_257_sink}
\vspace{-.1 cm}
\def\arraystretch{1.4}
\begin{tabular}{|c|c|c|c|c|c|c|c|c|c|c|c|c|c|c|c|c|c|c|c|c|c|c|}
\hline
$x_\mathrm{wall}$ & $0$ & $1$ & $2$ & $3$ & $4$ & $5$ & $6$ & $7$ & $8$ & $9$ & $10$ & $11$ & $12$ & $13$ & $14$ & $15$ & $16$ & $17$ & $18$ & $19$ & $20$ & $21$ \\ \hline
$v(x)$ & $0\atop 2$ & $0\atop 2$ & $1$ & $1$ & $0\atop 2$ & $2\atop 0$ & $1$ & $3$ & $2\atop 0$ & $2\atop 0$ & $0\atop 2$ & $3$ & $1$ & $0\atop 2$ & $0\atop 2$ & $1$ & $1$ & $2\atop 0$ & $2\atop 0$ & $3$ & $3$ & $2\atop 0$  \\[0.5mm] \hline
$x_\mathrm{sink}$ & $17$ & $16$ & $15$ & $14$ & $13$ & $12$ & $11$ & $10$ & $9$ & $8$ & $7$ & $6$ & $5$ & $4$ & $3$ & $2$ & $1$ & $22$ & $21$ & $20$ & $19$ & $18$ \\ \hline
\end{tabular}
\end{table}
\end{example}

Let us remark that the well-known Ferguson pairing principle \cite{berlekamp2004winning}
fails to hold in {\sc sink subtraction}.

\begin{rem}
Let us first reiterate the classical induction proof that, in the wall convention,
$v(x)=0$ if and only if $v(x+m)=1$, where $m=\min S$.
We then explain how this can fail in the sink convention.

Suppose $v(x)=0$ but $v(x+m)>1$. Then there exists a move $s$ such that
$v(x+m-s)=1$. Since $v(x)=0$, we have $v(x-s)>0$ for all wall-feasible $s$.
Hence subtracting $m$ from a heap of size $x+m-s$ does not yield a $0$-position,
which contradicts the induction hypothesis that every value~``1'' on a smaller heap
must be $m$-preceded by a~``0''.

Conversely, suppose $v(x+m)=1$ but $v(x)>0$, which implies $v(x)>1$.
Then there exists $s$ such that $v(x-s)=0$.
By induction, it follows that $v(x+m-s)=1$,
which contradicts the assumption that $v(x+m)=1$.

In the sink convention, the base case fails; see for example
Example~\ref{ex:257}, due to the presence of additional terminal zeroes.
Observe that $v(5)=1$, but $v(3)\ne 0$.
Indeed, the step ``By induction, then $v(x+m-s)=1$'' does not apply,
since in this case $x+m-s=0$, which belongs to the sink.
Similarly, $v(8)=0$ but $v(10)\ne 1$.
Here $v(8)$ is the first non-trivial zero-value,
yet $v(10)\ne 1$ because $v(5)=1$.
These violations of Ferguson's principle propagate to larger heaps.
Thus, in the sink convention, one can no longer deduce the
$1$-values from the $0$-values, or vice versa.
\end{rem}

\iffalse
Let us remark that the well known Ferguson's pairing principle \cite{berlekamp2004winning} fails to hold in {\sc sink subtraction}.

\begin{rem}
Let us first reiterate the classical induction proof that, in the wall convention, $v(x)=0$ if and only if $v(x+m)=1$, where $m = \min S$. Then we will see how this can fail in the sink convention. 

Suppose $v(x)=0$ but $v(x+m) > 1$. Then there is a move $s$ with $v(x+m-s)=1$. But, for all wall-feasible $s$, $v(x-s)>0$, so removing $m$ from a heap os size $x+m-s$ does not give $0$, which contradicts the induction hypothesis that every value ``1'' on a smaller heap must be $m$-preceeded by a ``0''. Next, suppose $v(x+m)=1$ but $v(x) > 0$, which implies $v(x)>1$. Then $v(x-s)=0$, for some $s$. By induction, then $v(x+m-s)= 1$, which contradicts the assumption. 

In the sink convention, the base case fails for example in Example~\ref{ex:257} due to the excess of terminal zeroes. Observe that $v(5)=1$, but $v(3)\ne 0$; indeed the part ``By induction, then $v(x+m-s)= 1$'' does not hold since in this case $x+m-s=0$, which belongs to the sink. Also $v(8)=0$ but $v(10)\ne 1$; in this case note that $v(8)$ is the first non-trivial zero-value, but $v(10)\ne 1$, since $v(5)=1$. The violations of Ferguson principle propagate up the heaps. So, in the sink convention, we may not be able to deduce the one-values from the zero-values or vice versa. 
\end{rem}
\fi
%Is there any other rule, similar to Ferguson's ``0-1 principle'', which is valid in the sink convention? 
Note that in the case of {\sc sink}, the above conjecture, together with Ferguson's pairing principle, instead implies a pairing of ``1''s and ``2''s: $v_\mathrm{ sink}(x)=1$ if and only if $v_\mathrm{sink}(x-m)=2$. 

\section{A proof of the main result}
\noindent The rest of the paper concerns the proof of Theorem~\ref{thm:main}. 
%\begin{proof}
We distinguish two cases:
\begin{enumerate}
    \item[(i)] $0 \le d \le m$;
    \item[(ii)] $m < d < 2m$.
\end{enumerate}
Note that, if $m=1$, only item~(i) applies. 
Let us prove {\bf item~(i)}. Write the moves as
\[
s_1 := m, \qquad s_2 := m+\delta, \qquad s_3 := 2m+\delta.
\]
%and let $d := d_2 = \delta \pmod{2m}$ with $0 \le d \le m$.

Define the {\em candidate period} as the word on the alphabet $\{0,1,2,3\}$: 
\begin{equation}\label{eq:easy_period}
(1^m2^m)^a\,3^{d}\,0^m\,(3^m0^m)^{a-1},
\end{equation}
where $a := \frac{\delta-d}{2m}+1$. The length of the candidate period \eqref{eq:easy_period} is
\begin{align*}
2ma+d+m+2m(a-1)
&=2m\!\left(\frac{\delta-d}{2m}+1\right)+d+m+2m\!\left(\frac{\delta-d}{2m}\right)\\
&=\delta-d+2m+d+m+\delta-d\\
&=3m+2\delta-d,
\end{align*}
as claimed.

Let us verify that the nim-sequence starting at $x=1$ follows the pattern of the candidate period 
\eqref{eq:easy_period}.

\smallskip
\noindent\emph{Case 1: the initial word $(1^m2^m)^a$.}
For $1 \le x \le 2ma = \delta-d+2m \le s_3$, we have $x-s_3 \le 0$, and hence, by the definition of the sink, 
$v(x-s_3)=0$. Since $s_1=m$, the mex rule forces alternating blocks of $1$ and $2$, each
of length $m$. The move $s_2=m+\delta$ does not interfere: the largest relevant options lie in
\[
[\,2ma-m-s_2,\; 2ma-s_2\,] = [-d,\; m-d] \subseteq (-\infty,m],
\]
whose nim-values are at most~$1$.
Thus $v(x)=2$ at the appropriate positions, completing the first part.

\smallskip
\noindent\emph{Case 2: the $3^d$ factor.}
For $2ma < x \le 2ma+d$, we have
\[
v(x-s_1)=2,\quad v(x-s_2)=1,\quad v(x-s_3)=0,
\]
since $x-s_3\le0$. Hence $v(x)=3$ for exactly $d$ consecutive positions.

\smallskip
\noindent\emph{Case 3: the $0^m$ factor.}
For $2ma+d < x \le 2ma+d+m = 3m+\delta$, all options satisfy $x-s_3>0$,
so no option has value~$0$. Thus the mex rule produces $m$ consecutive zeros.

\smallskip
\noindent\emph{Case 4: the trailing $(3^m0^m)^{a-1}$.} For the remaining $2m(a-1)$ positions, the same mechanism repeats.
The moves $s_3$ and $s_2$ together generate options of values $1$ and $2$, while $s_1=m$ toggles between producing and avoiding the value~$0$. Since $s_3-s_2=m$, it suffices to note that for all such $x$, both $x-s_3$ and $x-s_2$ lie in $[1,2ma]$, where the values are already known. This forces alternating blocks of $3^m$ and $0^m$.

\smallskip
Finally, observe that the next $2m+\delta$ positions reproduce the initial segment $(1^m2^m)^a3^d$. Hence the nim-sequence is purely periodic with period length $p(m,\delta)=3m+2\delta-d$. 
%\end{proof}

\begin{figure}[ht!]
     \centering
 \includegraphics[width=1\linewidth]{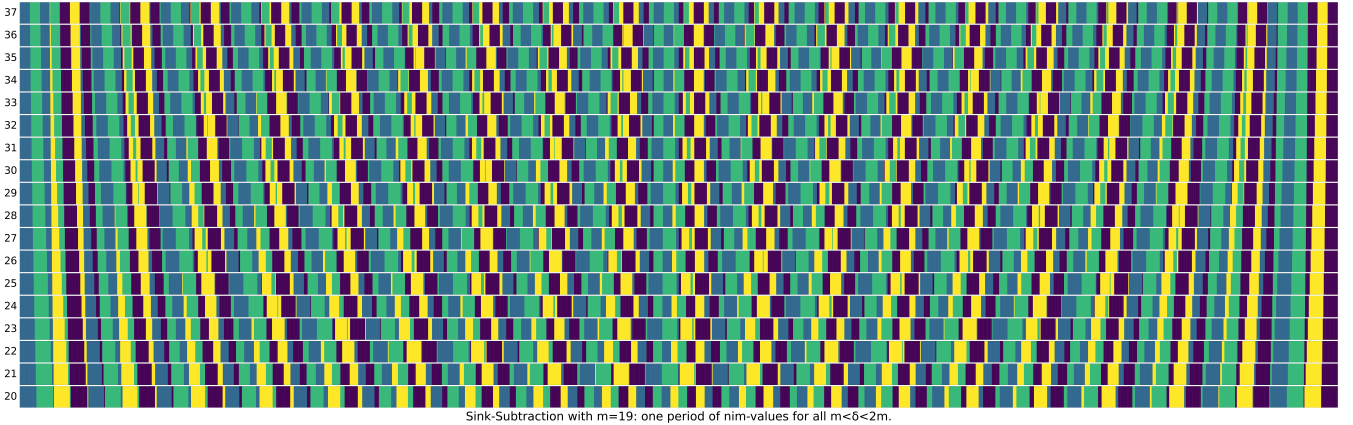}
     \caption{Nimber periods for $m=19$ and $0<k <19$, starting from $x=1$. The colors correspond to the nimbers; purple (0), blue (1), green (2) and yellow (3). The rows have different scaling as to fit the full periods in each case. Zoom in for details.}
     \label{fig:m19}
 \end{figure}

\medskip
\noindent Consider {\bf item (ii)}. For convenience we restrict $m< \delta <2m$, and let $\delta=m+k$. (The general case turns out to be an immediate consequence.) Similar to case (i), let the alphabet be $\{0,1,2,3\}$. As an example, let us begin by describing the periodic word for any rulesets of the form ``the lower case'', i.e. whenever $k=1$. The word is described by consecutive ``$B$-blocks'' of the form $$1^{m-i}0^i2^{m-i}1^{1+i}3^{m-1-i}2^{1+i}0^{m-i}3^{1+i},$$ for $0\le i\le m-1$, and at last, we concatenate the factor $0^m$.\footnote{The  period length for this case is mentioned in \cite{althofer1995superlinear}, but without giving the  explicit block construction or a proof.} By increasing the `block speed', the factor $1\le k<m$, another type of blocks appears, the ``$C$-blocks'', and the extreme ``upper case'' will favor only the $C$-blocks. The challenge is to understand the intermediate cases, with mixing between the block types; in Figure~\ref{fig:m19}, we illustrate all periods for the case of $m=19$, $0<k<19$.

Let $[x]_y$ denote the smallest nonnegative representative of $x\pmod y$. Let $[x]$ denote the smallest nonnegative representative of $x\pmod m$. Let $[x]^+$ denote the smallest positive representative of $x\pmod m$. For given $m, k$ and $i$, let $\alpha(i)=[ki]$, let $\beta(i)=[k(i+1)]$ and let $\gamma(i)=[k(i+1)]^+$. We drop $i$ in these notations, when the context allows.  

For given $1\le k<m$, $\delta=m+k$, for $0\le i\le m-1$, define the prefix word  
\begin{align}\label{eq:Ai}
A(i)=1^{m-\alpha}0^{\alpha}2^{m-\alpha}.
\end{align}

If 
%\begin{align}\label{eq:ki}
%    [k(i+1)]^+ > [ki].
%\end{align}
\begin{align}\label{eq:ki}
    \gamma > \alpha,
\end{align}
define the $i$th {\em block} as
%\begin{align}\label{eq:Bi}
%B(k,i)=x1^{[k(i+1)]^+}3^{m-[k(i+1)]^+}2^{[k(i+1)]^+}0^{m-[ki]}3^{[k(i+1)]^+},
%\end{align}
\begin{align}\label{eq:Bi}
B(i)=A(i)1^{\gamma}3^{m-\gamma}2^{\gamma}0^{m-\alpha}3^{\gamma},
\end{align}
and otherwise, if \eqref{eq:ki} does not hold, define the $i$th {\em block} as
%\begin{align}\label{eq:Ii}
%I(k,i)=x1^{m}2^{m}1^{[k(i+1)]}0^{k-[k(i+1)]}3^{m-k}2^{[k(i+1)]}3^{k-[k(i+1)]}0^m3^{[k(i+1)]}. 
%\end{align}
\begin{align}\label{eq:Ci}
C(i)=A(i)1^{m}2^{m}1^{\beta}0^{k-\beta}3^{m-k}2^{\beta}3^{k-\beta}0^m3^{\beta}. 
\end{align}

Observe that the $C$-block in \eqref{eq:Ci} has two exponents of the form $k-\beta$. Due to the reverse of inequality \eqref{eq:ki}, for a given $i$, $\alpha>\gamma\ge \beta$, then, by modular arithmetic, we must have $k>\beta$, so this exponent is always positive. 

For given $m$ and $k$, define the bijection $\varphi: [0,m-1]\rightarrow B^*$ induced by our construction, where $B^*$ is the set of all $B$- and $C$-blocks.
Let 
\begin{align}\label{eq:P}
P(m,k)=\prod_{0\le i< m} \varphi(i)\, \zeta(i), 
\end{align}
where $\zeta(i)=0^m$ if $\beta(i) =0$, and otherwise $\zeta(i)=\varnothing$. Note that $\beta =0$ whenever $i=m-1$, so $P$ always ends in the factor $0^m$. Indeed, if $m$ is prime, then $\zeta(i)=\varnothing$ if and only if $i<m-1$. Otherwise, we have additional inserted factors $0^m$, which are always preceded by a $B$-block, since $\beta=0$ implies $\gamma = m>\alpha$.

We verify that the nimbers in the additive subtraction game $S(m,2m+k, 3m+k)$ are purely periodic and the (not necessarily the shortest) period is $P = P(m,k)$, if $1<k<m$. That is, if $W=(w_i)_{i\ge 1}$ is the infinite word on the same alphabet such that, for all $i\ge 1$, $w_i$ is the nimber of $S$, then $W =PP\ldots$.

\begin{observation}
In spite of the periodicity hinging on a quadratic expression in $m$, we are able to verify correctness in constant time. Since each block consists of a constant number of factors, 8 and 12 respectively, and the number of interactions between blocks is also a constant, the task of testing the two mex properties (anti-collision and reachability) reduces to a case by case analysis. The number of tests is still large and a computer assisted proof seems like a convenient way forward. We have instead chosen to expose the data in a number of tables, so that every entry can be verified by a human. We explain by some examples how to interpret the tables. 
\end{observation}

The relevant differences $\gamma-\alpha$ and $\beta-\alpha$ do not depend on the index $i$, and consecutive blocks trigger some convenient identities. 
\begin{lemma}\label{lem:algebra}
    For given $i$, 
    \begin{enumerate}
        \item if $\gamma>\alpha$ then $\gamma-\alpha =k>0$, and otherwise,
        \item  if $\gamma<\alpha$, then $\alpha -\beta = m-k$, $\beta<k$ and $\alpha>m-k$. 
    \end{enumerate}
        Moreover, for $i\in \{0,\ldots ,m-1\}, \alpha(i)=\beta(i-1)$ and for $i\in \{1,\ldots ,m\}, \alpha(i)=\gamma(i-1)$. 
\end{lemma}
\begin{proof}
    Apply modular arithmetic.
\end{proof}
This lemma will be used often without explicit mention, as we regard these relations as basic facts that control much of the mex behavior of the rulesets. %In the second item, we will have use for $\beta$ 

\begin{figure}[ht]
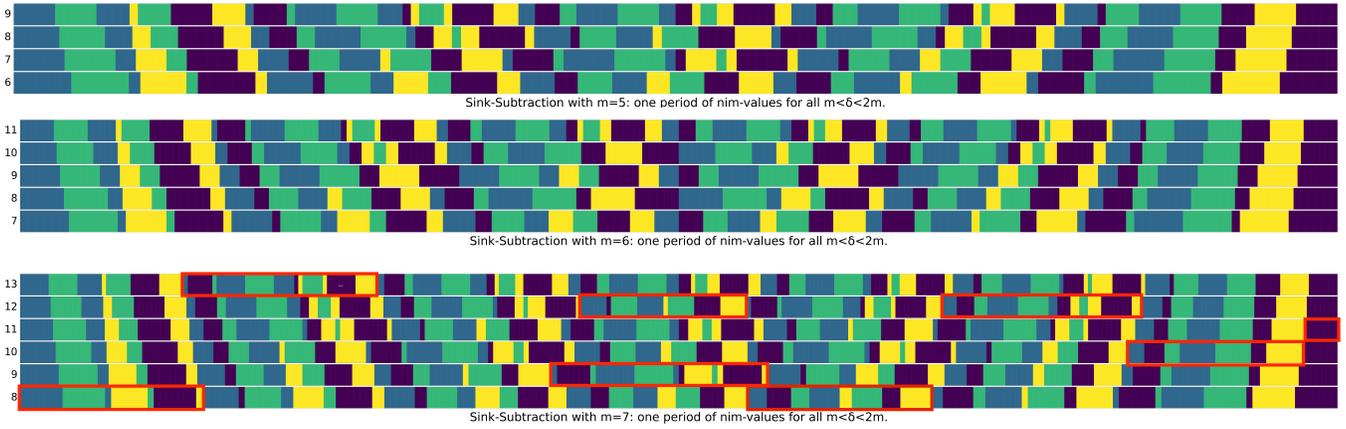

    \centering
    \includegraphics[width=1\linewidth]{m5allnosink.png}\vspace{-.2cm}
    \includegraphics[width=1\linewidth]{m6allnosink.png}
    \includegraphics[width=1\linewidth]{m7all_nosink_red.png}
    \caption{Nimber periods for $m=5,6,7$ and $0<k <m$, respectively, starting at $x=1$. The rows have different scaling as to fit the full periods in each case. In the case of $m=6$, and rows $k=2,3,4$ the smallest period is a divisor of the full row length. Hence in these cases there are more than one $\zeta$-factors. In the case of $m=7$, we have indicated with red rectangles, for $k=1$ the blocks $B(0)$ and $B(4)$; for $k=2$ the block $C(3)$; for $k=3$ the block $B(6)$; for $k=4$ the factor $\zeta(6)$; for $k=5$, the blocks $B(3)$ and $C(5)$; for $k=6$ the block $C(1)$.}
    \label{fig:m5m6}
\end{figure}

Let us prove that the construction produces the correct period lengths. For example, if $m=5$ and $k=1$, then (by ignoring technical details) the block structure is $BBBBB0^m$, which satisfies the claimed period length $5(4\cdot 5+2c-2a)+5=115$. If $m=5$ and $k=4$, then the block structure is $BCCCB0^m$, which satisfies the claimed period length $2(4\cdot 5+8)+3(6\cdot 5+4-b+a)+5=56+99+5=160$. The cases for $m=5$ are illustrated in the upper picture of Figure~\ref{fig:m5m6}. If $m=6$ and $k=2$, then the block structure is $BBB0^mBBB0^m$, with claimed period length $3(4\cdot 6+2c-2a)+6=90$. If $m=6$ and $k=4$, then the block structure is $BCB0^mBCB0^m$, with claimed period length $2(4\cdot 6+2c-2a)+6\cdot6+4+b-a+6=108$. The cases for $m=6$ are illustrated in the middle picture of Figure~\ref{fig:m5m6}. Recall that we here assume $\delta = m+k$, $1\le k<m$. 

%%%%%%%%%%Figures m=5,6 here

%%%%%%%%%%%%%%%%%%%%%%%%%%
\begin{lemma}
    For given $m,k$, we have $|P(m,k)| = m(4m+3k)/\gcd(m,k)$.
\end{lemma}
\begin{proof}
Let $g := \gcd(m,k)$. Write $m = gM$ and $k = gK$ with $\gcd(M,K)=1$. Recall the definitions
$\alpha(i) = [ki]_m, \beta(i) = [k(i+1)]_m$ and  $\gamma(i) = [k(i+1)]_m^+$. Moreover, for each $i$, we build a block $B(i)$ or $C(i)$ with prefix $A(i)$, and appending an extra factor $\zeta(i)=0^m$ to a $B$-block if $\beta(i)=0$.

\medskip\noindent
\textbf{1)} Periodicity in the block index $i$. For $i\ge 0$,
\begin{align*}
\alpha(i+M) &= [k(i+M)]_m\\ &= [ki + kM]_m\\ &= [ki + (km/g)]_m\\ &= [ki]_m\\ &= \alpha(i),
\end{align*}
since $g|k$.
The same calculation shows $\beta(i+M)=\beta(i)$ and hence also $\gamma(i+M)=\gamma(i)$.

Thus the whole pattern of blocks (including the presence or absence of
$\zeta(i)=0^m$) is periodic in $i$ with period $M$.  In particular, the word
$P(m,k)$ is obtained by concatenating the $M$ blocks
corresponding to $i=0,\dots,M-1$, repeated $g$ times.  Since we are interested
in the length of one period, it suffices to count the total length contributed
by a single index cycle $i=0,\dots,M-1$.

\medskip\noindent
\textbf{2)} Block lengths. From \eqref{eq:Ai}–\eqref{eq:Ci}, the prefix $A(i)=1^{m-\alpha}0^{\alpha}2^{m-\alpha}$ has length $|A(i)|=2m-\alpha(i)$.

A direct inspection of the definitions of $B(i)$ and $C(i)$ yields:
\begin{enumerate}
    \item If $\gamma(i) > \alpha(i)$ and $\beta(i)>0$ (the generic $B$-case),
    then $B(i)$ has length $|B(i)| = 4m+2k.$
    \item If $\gamma(i) > \alpha(i)$ and $\beta(i)=0$ (the “wrap” $B$-case
    with an inserted $0^m$), or if $\gamma(i) < \alpha(i)$ (the $C$-case),
    then the resulting `block' (either $B(i)\zeta(i)$ or $C(i)$) has length $5m+2k.$
\end{enumerate}
So every index $i$ in one cycle contributes either a \emph{short block} of
length $4m+2k$ or a \emph{long block} of length $5m+2k$.

\medskip\noindent
\textbf{3)} Classification via $\alpha(i)$. We claim that the type of block is determined solely by whether
$\alpha(i) < m-k$ or $\alpha(i) \ge m-k$. 

Indeed, if $\alpha(i) < m-k$, then $\alpha(i) + k < m$, 
so $k(i+1) = ki + k$ does not wrap modulo $m$, and thus 
\[
\beta(i) = \gamma(i) = \alpha(i)+k > \alpha(i).
\]
Hence we are in the generic $B$-case and obtain a short block of length
$4m+2k$. 

Conversely, if $\alpha(i) \ge m-k$, then $\alpha(i)+k \ge m$, and so
%\begin{align*}
%k(i+1) = ki + k \equiv \alpha(i)+k-m \pmod m
%\end{align*}
$\beta(i)$ wraps around $m$, which implies $\beta(i) = \alpha(i)+k-m < \alpha(i)$. 
If $\beta(i)=0$, we are in the wrapped $B$-case with an extra
$\zeta(i)=0^m$; if $\beta(i)>0$, then $\gamma(i)=\beta(i)<\alpha(i)$ and we
are in the $C$-case.  In both subcases we obtain a  long block of length
$5m+2k$. Thus:

\[
\begin{cases}
\alpha(i) < m-k &\Rightarrow\ \text{short block of length }4m+2k,\\[2pt]
\alpha(i) \ge m-k &\Rightarrow\ \text{long block of length }5m+2k.
\end{cases}
\]

\medskip\noindent
\textbf{4)} The number of short/long blocks in one cycle. 
Over one index cycle $i=0,\dots,M-1$, the values of $\alpha(i)$ run
through the multiples of $g$ modulo $m$:
$$\{\alpha(i)\mid 0\le i\le M-1\} = \{0,g,2g,\dots,m-g\},$$ 
each exactly once. Write $m-k = g(M-K)$. 
Then the inequality $\alpha(i) \ge m-k$ is equivalent to
\[
gj \ge g(M-K)
\quad\Longleftrightarrow\quad
j \ge M-K,
\]
where $\alpha(i) = gj$ with $0\le j\le M-1$.
Hence, in one cycle:
\begin{enumerate}
\item there are exactly $M-K$ indices with $\alpha(i) < m-k$ (short blocks), and  
\item there are exactly $K$ indices with $\alpha(i) \ge m-k$ (long blocks).
\end{enumerate}
\medskip\noindent
\textbf{5)} The length of $P$. Since $P$ consists exclusively of short and long blocks, we have that 
%Let $S := M-K$ and $L:=K$ be the numbers of short and long blocks respectively in a single cycle. Then
\begin{align*}
|P(m,k)| &= (M-K)(4m+2k) + K(5m+2k)\\
&= M(4m+2k) + Km.
\end{align*}
By substituting $M=m/g$ and $K=k/g$ we get 
\begin{align*}
|P(m,k)|
&= \frac{m(4m+2k)}{g} + \frac{km}{g}\\
&= \frac{m(4m+3k)}{g}\\
&= \frac{m(4m+3k)}{\gcd(m,k)}.
\end{align*}
\end{proof}

%\begin{lemma}
%    For given $m,k$, then $|P(m,k)|=m(4m+3k)/\gcd(m-k,k)$. 
%\end{lemma}
%\begin{proof}
    
%\end{proof}
%We can express as pseudocode as follows.

%\paragraph{Indexing convention.}

We use elementary algebra to justify mex-consistency.  Each {\em maximum factor} $F$, of $P$, is of the form $F=a^b$ with $a\in\{0,1,2,3\}$ and where $0\le b\le m$ is maximum such that adjacent letters differ from $a$; the {\em length} of $F$ is $|F| = b$. We often omit the word `maximum' and call a maximum factor simply by a {\em factor}. The {\em factor-value} of a factor $F$ is $v(F)=a$ if $F=a^b$. 
We will treat each factor as a single entity, and use terminology such as {\em factor-option} and {\em factor-value}. The advantage is that it drastically reduces the number of cases, but an apparent disadvantage is that sometimes individual positions within a factor can have options in distinct factors, for a given move. This phenomenon proves unavoidable though, and is part and parcel of the quadratic period length, as we will see. By intuition, when factor relations settle, then we approach the full period. 

The blocks will use two independent indices, the latter of which will enumerate its factors.
\begin{enumerate}
    \item The \emph{block index} $i$ refers to the position of a block in the
    concatenated word $P(m,k)$. Thus $B(i)$ and $C(i)$ denote the $i$th block,
    of type $B$ or $C$, respectively.

    \item Within a fixed block, we label its consecutive \emph{factors} by a
    local index. For a given $B(i)$-block, the factors are denoted
    $A_1,A_2,A_3, B_1,\ldots,B_5$, where $A_j$ denotes the $j$th factor of the prefix $A(i)$ and $B_j$ refers to the
    $(j+3)$rd factor of that block.
\end{enumerate}
Unless explicitly stated otherwise, the block index $i$ is fixed when referring to factor indices. For example, the notation $B_3$ always refers to the sixth factor within a single $B(i)$ block, where the block index $i$ is given by the surrounding context. Whenever we refer to an object other than a block or a factor, we explicitly use the term \emph{word}.

%A note on indexing of blocks. While we study the blocks there is a double indexing. The $i$-index counts the current block's position in the word $P(m,k)$, $B(i)$ or $C(i)$, while the $j$ in \texttt{Bj} refers to the $j+3$rd factor of the given $B$ block, where $j\in \{1,\ldots, 5 \}$. Similarly \texttt{Aj} denotes the $j$th factor of the $A(i)$ prefix (for a given $i$). 

Recall that the block length of $B$ is $4m-2\alpha +2\gamma =4m+2k>3m+k$ (by $\gamma>\alpha$), while the block length of $C$ is $6m+k-\alpha+\beta=5m+2k>3m+k$ (by $\beta<\alpha$). 
One can check that, under block concatenation, the distance between two {\em similar factors} is always greater than the largest move, $3m+k$; two $A$-factors are {\em similar}, if they carry the same $j$-index, and analogously for $B$- and $C$-factors. 

We display tables with factor-options and their values. Each table corresponds to a local block concatenation and verifies one of the two mex conditions. Tables~\ref{tab:sinkB} to \ref{tab:CC} justify {\em reachability}, that all smaller values appear as values of options of the conjectured factor-values. Here, of course, we use standard induction principles. Tables~\ref{tab:ACsinkB} to \ref{tab:antiCC} verify {\em anti-collision}, the second mex-principle, that any conjectured factor-value, cannot have the same value at the level of any factor-option. 

%In the tables, rather than fixing a specific heap of size one, we adopt the following convention. For each row (i.e.\ each game value under consideration), we choose as reference point the smallest position in the leftmost factor of that value; for the sink we use position $0$ as reference. 

We often study the corresponding set of positions of a given factor. To this purpose it is convenient to temporarily identify {\em position one} with the first letter of some relevant factor, and we use a special font \texttt{F} to denote the {\em factor set}, the set of such induced positions corresponding to a factor $F$. Any such notation is local and used as a computational device. Each factor is then interpreted as a set of consecutive positions, all carrying the same game value, determined relative to this reference point. In this way, the symbolic word representations, together with the specified factor lengths, are mapped unambiguously to actual game positions, allowing us to analyze move options between factors.

Throughout the tables, the symbols $s_1$, $s_2$, and $s_3$ denote the moves $m$, $2m+k$, and $3m+k$, respectively. In all reachability tables, a dash “--” indicates that no reachability condition is required from the corresponding factor to that value. In tables for concatenations (e.g.\ \texttt{bB} concerning reachability in \eqref{eq:Bi}), we use an uppercase letter for the current block $B(i)$ and a lowercase for the preceding block $B(i-1)$. For example, the factor set \texttt{B4} is the set of positions corresponding to the factor $0^{m-\alpha}$ in $B(i)$, while \texttt{b4} is the set of positions corresponding to the factor $0^{m-\alpha'}$ in $B(i-1)$; throughout the tables, we write $\alpha'=\alpha(i-1)$, $\beta'=\beta(i-1)$ and $\gamma'=\gamma(i-1)$. 
In the reachability tables, distances are measured from the beginning of the
relevant factor to the end of the preceding block.  For example,  the distance to the end of the $B$-block of the factor
\texttt{b4} is 
$d(\texttt{b4}) = m-\gamma' + \alpha' = m-k$ (by Lemma~\ref{lem:algebra}),  
and its length is $|\texttt{b4}|=\gamma(i-1)=\gamma'$. Similarly, $d(\texttt{c4}) = 2m+k
\text{ and } 
d(\texttt{c8}) = m+\beta', 
$
with corresponding lengths
$
|\texttt{c4}| = k-\beta'
\text{ and } 
|\texttt{c8}| = m.
$

%Throughout, ``factor'' refers to a symbolic word of the form $x^y$, while ``factor set'' refers to the induced set of game positions when a reference point is fixed. We only pass to factor sets when verifying reachability or anti-collision. 
The value of a factor set is its factor-value. A {\em factor-option} of $F$ is a factor $F'$ such that there is an $x\in \texttt{F'}$ such that, for some $y\in \texttt{F}$, $x-y\in \{m, 2m+k, 3m+k\}=\{s_1,s_2,s_3\}$. To read the tables: fix a column factor set. Each row lists all potential factor-options of a relevant game value.  Reachability tables ensure the existence of smaller values as move options, while anti-collision tables exclude equal values in the set of move options; here a triple of inequalities corresponds to the moves $s_1, s_2, s_3$, respectively.

%The verification of the periodic pattern proceeds by a finite case analysis over possible block concatenations. Specifically, we consider the concatenations \texttt{sinkB}, \texttt{bZB}, \texttt{bB}, \texttt{bC}, \texttt{cB}, and \texttt{cC}. For each case we verify two conditions: \emph{reachability}, ensuring that each position has an option to all smaller game values, and \emph{anti-collision}, ensuring that no move maps a position to another position of the same value. The reachability conditions are summarized in Tables~\ref{tab:sinkB}--\ref{tab:CC}, and the anti-collision conditions in Tables~\ref{tab:ACsinkB}--\ref{tab:antiCC}. Throughout these verifications, the algebraic relations collected in Lemma~\ref{lem:algebra} are used repeatedly to reduce distance comparisons.

%By using Lemma~\ref{lem:algebra}

\medskip 
\noindent {\bf Example \texttt{A3} in \texttt{sinkB} vs. \texttt{bZB}.} This is the case where $\beta(i-1)=\alpha(i)=0$, and hence $\zeta(i-1)=0^m$. Moreover, it follows that $\alpha(i-1)=m-k$. The factor set \texttt{A3} has size $m$ in such $B(i)$, and, by \eqref{eq:Bi}, we study the word 
\begin{align}
B(i-1)\zeta(i-1)B(i)=1^{k}0^{m-k}2^{k}1^{m}3^{0}2^{m}0^{k}3^{m}\underline{0}^m1^{m}0^{0}\boldsymbol{2}^{m}1^{k}3^{m-k}2^{k}0^{m}3^{k}, 
\end{align}
where the underlined 0-factor is $\zeta(i-1)$. Let us verify the factor set \texttt{A3}, induced by the bold font factor $\boldsymbol{2}^{m}$, assuming by induction that the previous factors are correctly defined and satisfy the mex-rule. 
 Specifically, we have to prove that, for all positions in \texttt{A3}, the value is $v(2^{m})=2$. 
 
 By inspection, and since $\texttt{A2}=\emptyset$, all the $m$-options have value one, by induction. Note that in the case of sink preceding $B(0)$ (denoted \texttt{sinkB}), the related justification of \texttt{A3} is immediate, since the move $m$ assures the value one, and both larger moves drown in the sink. So, let us focus on the ``wrap case'' of $B(i-1)\zeta (i-1)B(i)$, also denoted \texttt{bZB}; see Tables~\ref{tab:BZB} and \ref{tab:antiBZB} for general reachability and anti-collision justifications, respectively. 

 What differs from the sink argument is that both long options have to be scrutinized for reachability (to find value 0) and for anti-collision (to avoid value 2). 

Let us denote by $\langle x \rangle =\{1,\ldots ,x\}$. The claim from Table~\ref{tab:BZB} (reachability) is that both long factor-options contribute to the value zero. In fact, \texttt{Z} alone cannot compensate for the loss of some sink options. Instead, we use Table~\ref{tab:BZB} to compute the `inverse images' of both 0-factor-options \texttt{Z} and \texttt{b4}, and prove that their union is a superset of \texttt{A3}. We may take as reference point the first zero-value in \texttt{b4} as position `one'. The 0-factor set \texttt{b4} has size $m - \alpha' = m-(m-k) = k$, while \texttt{b5} has size $\gamma = [km]^+=m$ and the 0-factor set \texttt{Z}, which starts at position $m+k$, has size $m$. Thus, the candidate super-set is 
\begin{align*}\label{eq:Gamma}
\Gamma&=(s_3+\langle k \rangle)\cup (m+k+s_2+\langle m \rangle) \\ 
&= (3m+k+\langle k \rangle)\cup (m+k+2m+k+\langle m \rangle) \\
&= \{3m+k+1,\ldots ,3m+2k,3m+2k+1,\ldots,4m+2k\}.
\end{align*}
By using the same reference point, the positions of \texttt{A3} are $\{3m+k+1,\ldots ,4m+k \}\subset \Gamma$; namely,  since $|\texttt{A1}|=m-\alpha = m$, the position just before \texttt{A3} is $k+m+m+m=3m+k$. We may conclude that all positions in \texttt{A3} have both values 0s and 1s as options. 

Next we verify that, as claimed by Table~\ref{tab:antiBZB} (anti-collision), the intersection of move-translated 2-factor positions with  \texttt{A3} is empty. 
%Since the block length of $B$ is greater than $3m+k$, positions in \texttt{a3} cannot be options from \texttt{A3}. 
There is exactly one 2-factor between \texttt{a3} and \texttt{A3}, namely \texttt{b3}. %which has length $\gamma' = [km]^+=m$. 
The number of letters between \texttt{b3} and \texttt{A3} is 
\begin{align}
|\texttt{b4}|+|\texttt{b5}|+|\texttt{Z}|+|\texttt{A1}|+|\texttt{A2}|&=
m-\alpha(m-1)+ \gamma(m-1)+\boldsymbol{m}+ m-\alpha(0)+ m-\gamma(0)\\
&=k+m+m+m+m-k\\
&=4m\\ 
&>3m+k\\
&=\max(S).
\end{align}
Here the bold $m$ counts the number of symbols from the inserted $\zeta(i-1)=0^m$ factor. A consequence is that \texttt{a3} is also sufficiently separated.  
By induction, this proves the correctness of the mex-values for the factor \texttt{A3} in case \texttt{bZB}.\\

\noindent {\bf Example \texttt{A3} in \texttt{bB}.} This is the general case of consecutive $B$-blocks; see Table~\ref{tab:BB}. We study the 2-factor set \texttt{A3} in the word $B(i-1)B(i)$, 
\begin{align}\label{eq:Bi}
 1^{\gamma'}3^{m-\gamma'}2^{\gamma'}0^{m-\alpha'}3^{\gamma'}\underline{1}^{m-\alpha}0^{\alpha}\boldsymbol{2}^{m-\alpha}1^{\gamma}3^{m-\gamma}2^{\gamma}0^{m-\alpha}3^{\gamma},
\end{align}
where we have omitted the prefix $A(i-1)$, since it is separated by a block from \texttt{A3}. Similar to the \texttt{bZB} case, the move $s_1$ provides all $m-\alpha$ $1$-values (note that in this case \texttt{A2} is non-empty). But here, only one move suffices to reach the $m-\alpha$ $0$-values, namely $s_2$.  The candidate superset is 
%\begin{align}\label{eq:bBA3}
$\Gamma = s_2 + \langle m -\alpha'\rangle$.
%\\ &= \{2m+k+1,\ldots , 3m+k-\alpha'\}.
%\end{align}
 By using the same reference position, the set \texttt{A3} is: 
 \begin{align}\label{eq:bBA3}
 m-\alpha'+\gamma' + m - \alpha + \alpha + \langle m - \alpha\rangle &= 2m+k+\langle m - \alpha\rangle\\
 &\subset\Gamma, 
 \end{align}
 by Lemma~\ref{lem:algebra}, and since $\alpha>\alpha'$ (recall this is a case where $\beta(i-1) = \alpha(i)>0$). For anti-collision, by \eqref{eq:Bi}, the only $2$-factor set we have to test is \texttt{b3}, of size $\gamma'$; see Table~\ref{tab:antiBB}. We show that \texttt{b3} fits between the moves $s_2$ and $s_3$. By Lemma~\ref{lem:algebra}, the number of letters  between \texttt{b3} and \texttt{A3} is $2m+k$, which proves the first part. But $\gamma' = \alpha$, so $\gamma'+2m+k+m-\alpha=3m+k$, which implies that \texttt{A3} does not collide with \texttt{b3} via the move $s_3$. 
 
 Observe that this example resolves also \texttt{A3} in \texttt{bC}, since $A(i)$ is a prefix in both the $B$ and $C$ blocks.\\

\noindent {\bf Example \texttt{A3} in \texttt{cC}.} Since $A(i)$ is a prefix in both the $B$ and $C$ blocks, this example resolves both \texttt{A3} in \texttt{cB} and \texttt{A3} in \texttt{cC}. Curiously, there are double entries for reachability in Table~\ref{tab:CC} (and the same in Table~\ref{tab:CB}). The claimed moves for reachability to the two 0-factor sets \texttt{c8} and \texttt{c4} are $s_2$ and $s_3$, respectively. There are $m+\beta'$ letters between \texttt{c8} and \texttt{A3}. The candidate superset is 
\begin{align}
\Gamma= (s_3+\langle k-\beta'\rangle)\cup (m+k-\beta'+s_2 +\langle m\rangle),
\end{align}
with $\beta'<k$, since we study a $C$-block. Moreover, since  $\alpha=\beta'$, we get the interval 
\begin{align}
\Gamma &= \{3m+k+1,\ldots , 3m+k+k-\alpha\}\cup\{3m+2k-\alpha+1,\ldots , 4m+2k-\alpha\}\\
&=\{3m+k+1,\ldots , 4m+2k-\alpha\}. 
\end{align}
And, with the same reference point, by inspection,
\begin{align}
2m+\beta' + m+k-\beta'+\texttt{A3}&=3m+k+\langle m-\alpha\rangle\\ &=\{3m+k+1,\ldots ,4m+k-\alpha\}\subset \Gamma.
\end{align}
Thus reachability holds, and note that both moves $s_3$ and $s_2$ are required to $0$-cover the factor. 

Next, for anti-collision, as displayed in Table~\ref{tab:antiCC}, the relevant 2-factors are  \texttt{c2} and \texttt{c6}. But, by the table,  \texttt{c2} is distanced by more than $3m+k$, while \texttt{c6} is distanced exactly by $2m+k$, and has size $m-\alpha < m$, and hence neither $s_3$ collides with \texttt{A3} in \texttt{cC}.\\

\noindent {\bf Example \texttt{C9}.} In this example, we guide the reader through table justification. Note that \texttt{C9} is distanced by more than $4m$ to its $A$-prefix. Hence, both reachability and anti-collision will be resolved locally within the sets \texttt{C1} to \texttt{C9}. Thus, here we resolve this entry simultaneously in all four tables for \texttt{bC} and \texttt{cC}. Indeed, all Tables~\ref{tab:BC},~\ref{tab:CC},~\ref{tab:antiBC} and~\ref{tab:antiCC} have the same \texttt{C9} entries. For reachability, $s1$ maps \texttt{C9} into \texttt{C8}, $s2$ maps \texttt{C9} onto \texttt{C3} and $s3$ maps \texttt{C9} into \texttt{C2}. The relevant sizes and distances are displayed in those tables, so verification is direct. Thus \texttt{C9} satisfies reachability. 

For anti-collision, both \texttt{C5} and \texttt{C7} are 3-factor sets, and require tests; see the last column and last two rows in Table~\ref{tab:BC}. Observe that the size of $|\texttt{C8}|=m$, so $s_1$ does not collide with \texttt{C7} and hence also not with \texttt{C5}. For $s_2$, note that $|\texttt{C5} \cup \texttt{C6} \cup \texttt{C7}\cup \texttt{C8} | = 2m$, and hence the two last inequalities are also correct in both cells (since $k\ge 1$). Thus \texttt{C9} does not collide. 
\\

\begin{table}[ht!]
\begin{center}
\caption{{\bf Reachability \texttt{sinkB}.} This is the base case, \texttt{sinkB}, where the block $B(0)$ meets the sink. We omit the \texttt{A2} column, since $\alpha(0)=0$. For each column factor, the table lists moves and target factors whose inverse
images jointly cover the entire factor, ensuring reachability of the indicated
value.} %The ``-'' squares do not require a reachability test.}
\label{tab:sinkB}
\vspace{-3pt}
\setlength{\tabcolsep}{3pt}
\renewcommand{\arraystretch}{1.15}
\begin{tabular}{|c|c|c|c|c|c|c|c|}
\hline
\texttt{sinkB}  & $\underset{m}{\texttt{A1} \text{``1''}}$ & $\underset{m}{\texttt{A3} \text{``2''}}$ & $\underset{k}{\texttt{B1} \text{``1''}}$ & $\underset{m-k}{\texttt{B2} \text{``3''}}$ & $\underset{k}{\texttt{B3} \text{``2''}}$ & $\underset{m}{\texttt{B4} \text{``0''}}$ & $\underset{k}{\texttt{B5} \text{``3''}}$ \\
\hline\hline
$ 0 $ & \begin{tabular}{@{}l@{}} $s_1$ \\ \texttt{sink} \end{tabular} & \begin{tabular}{@{}l@{}} $s_2$ \\ \texttt{sink} \end{tabular} & \begin{tabular}{@{}l@{}} $s_2$ \\ \texttt{sink} \end{tabular} & \begin{tabular}{@{}l@{}} $s_3$ \\ \texttt{sink} \end{tabular} & \begin{tabular}{@{}l@{}} $s_3$ \\ \texttt{sink} \end{tabular} & \raisebox{1.5ex}{--} & \begin{tabular}{@{}l@{}} $s_1$ \\ \texttt{B4} \end{tabular} \\
\hline
$ 1 $ & \raisebox{1.5ex}{--} & \begin{tabular}{@{}l@{}} $s_1$ \\ \texttt{A1} \end{tabular} & \raisebox{1.5ex}{--} & \begin{tabular}{@{}l@{}} $s_2$ \\ \texttt{A1} \end{tabular} & \begin{tabular}{@{}l@{}} $s_1$, $s_2$ \\ \texttt{B1},\ \texttt{A1} \end{tabular} & \raisebox{1.5ex}{--} & \begin{tabular}{@{}l@{}} $s_2$ \\ \texttt{B1} \end{tabular} \\
\hline
$ 2 $ & \raisebox{1.5ex}{--} & \raisebox{1.5ex}{--} & \raisebox{1.5ex}{--} & \begin{tabular}{@{}l@{}} $s_1$ \\ \texttt{A3} \end{tabular} & \raisebox{1.5ex}{--} & \raisebox{1.5ex}{--} & \begin{tabular}{@{}l@{}} $s_3$ \\ \texttt{A3} \end{tabular} \\
\hline
\end{tabular}
\end{center}
\end{table}

% === BZB ===
\begin{table}[ht!]
\begin{center}
\caption{{\bf Reachability \texttt{bZB}.}  This is the wrap case, where the block $B(i)$ meets the inserted 0-factor $\zeta(i-1)=0^m$ preceded by $B(i-1)$. We omit the \texttt{A2} column, since $\alpha(0)=0$. For each column factor, the table lists moves and target factors whose inverse images jointly cover the entire factor, ensuring reachability of the indicated value. For the factor \texttt{b4}, the distance to the current $B$-block is $d(\texttt{b4}) = \gamma(m-1) = m$, and its length is $|\texttt{b4}|=m-\alpha(m-1)=k$.
}
\label{tab:BZB}
\vspace{-3pt}
\setlength{\tabcolsep}{3pt}
\renewcommand{\arraystretch}{1.15}
\begin{tabular}{|c|c|c|c|c|c|c|c|}
\hline
\texttt{bZB}
& $\underset{m}{\texttt{A1} \text{``1''}}$
%& $\underset{0}{\texttt{A2} \text{``0''}}$
& $\underset{m}{\texttt{A3} \text{``2''}}$
& $\underset{k}{\texttt{B1} \text{``1''}}$
& $\underset{m-k}{\texttt{B2} \text{``3''}}$
& $\underset{k}{\texttt{B3} \text{``2''}}$
& $\underset{m}{\texttt{B4} \text{``0''}}$
& $\underset{k}{\texttt{B5} \text{``3''}}$ \\
\hline\hline
$0$
& \begin{tabular}{@{}l@{}} $s_1$ \\ \texttt{Z} \end{tabular}
%& \raisebox{1.5ex}{--}
& \begin{tabular}{@{}l@{}} $s_2$, $s_3$ \\ \texttt{Z},\ \texttt{b4} \end{tabular}
& \begin{tabular}{@{}l@{}}  $s_2$ \\ \texttt{Z} \end{tabular}
& \begin{tabular}{@{}l@{}} $s_3$ \\ \texttt{Z} \end{tabular}
& \begin{tabular}{@{}l@{}}  $s_3$ \\  \texttt{Z} \end{tabular}
& \raisebox{1.5ex}{--}
& \begin{tabular}{@{}l@{}} $s_1$ \\ \texttt{B4} \end{tabular} \\
\hline
$1$
& \raisebox{1.5ex}{--}
%& \raisebox{1.5ex}{--}
& \begin{tabular}{@{}l@{}} $s_1$ \\ \texttt{A1} \end{tabular}
& \raisebox{1.5ex}{--}
& \begin{tabular}{@{}l@{}} $s_2$ \\ \texttt{A1} \end{tabular}
& \begin{tabular}{@{}l@{}} $s_1$, $s_2$ \\ \texttt{B1},\ \texttt{A1} \end{tabular}
& \raisebox{1.5ex}{--}
& \begin{tabular}{@{}l@{}} $s_2$ \\ \texttt{B1} \end{tabular} \\
\hline
$2$
& \raisebox{1.5ex}{--}
%& \raisebox{1.5ex}{--}
& \raisebox{1.5ex}{--}
& \raisebox{1.5ex}{--}
& \begin{tabular}{@{}l@{}} $s_1$ \\ \texttt{A3} \end{tabular}
& \raisebox{1.5ex}{--}
& \raisebox{1.5ex}{--}
& \begin{tabular}{@{}l@{}}  $s_3$ \\  \texttt{A3} \end{tabular} \\
\hline
\end{tabular}
\end{center}
\end{table}

% === BB (modern/base style; correct undersets) ===
\begin{table}[ht!]
\begin{center}
\caption{{\bf Reachability \texttt{bB}.} This is the generic case where the block $B(i)$ meets $B(i-1)$. For the factor \texttt{b4}, the distance to the current $B$-block is $d(\texttt{b4}) = \gamma'$, and its length is $|\texttt{b4}|=m-\alpha'$.}
\label{tab:BB}
\vspace{-3pt}
\setlength{\tabcolsep}{3pt}
\renewcommand{\arraystretch}{1.15}
\begin{tabular}{|c|c|c|c|c|c|c|c|c|}
\hline
\texttt{bB}
& $\underset{m-\alpha}{\texttt{A1} \text{``1''}}$
& $\underset{\alpha}{\texttt{A2} \text{``0''}}$
& $\underset{m-\alpha}{\texttt{A3} \text{``2''}}$
& $\underset{\gamma}{\texttt{B1} \text{``1''}}$
& $\underset{m-\gamma}{\texttt{B2} \text{``3''}}$
& $\underset{\gamma}{\texttt{B3} \text{``2''}}$
& $\underset{m-\alpha}{\texttt{B4} \text{``0''}}$
& $\underset{\gamma}{\texttt{B5} \text{``3''}}$ \\
\hline\hline
$0$
& \begin{tabular}{@{}l@{}} $s_1$ \\ \texttt{b4} \end{tabular}
& \raisebox{1.5ex}{--}
& \begin{tabular}{@{}l@{}} $s_2$ \\ \texttt{b4} \end{tabular}
& \begin{tabular}{@{}l@{}} $s_1$, $s_2$ \\ \texttt{A2},\ \texttt{b4} \end{tabular}
& \begin{tabular}{@{}l@{}} $s_3$ \\ \texttt{b4} \end{tabular}
& \begin{tabular}{@{}l@{}} $s_2$, $s_3$ \\ \texttt{A2},\ \texttt{b4} \end{tabular}
& \raisebox{1.5ex}{--}
& \begin{tabular}{@{}l@{}} $s_1$, $s_3$ \\ \texttt{B4},\ \texttt{A2} \end{tabular} \\
\hline
$1$
& \raisebox{1.5ex}{--}
& \raisebox{1.5ex}{--}
& \begin{tabular}{@{}l@{}} $s_1$ \\ \texttt{A1} \end{tabular}
& \raisebox{1.5ex}{--}
& \begin{tabular}{@{}l@{}} $s_2$ \\ \texttt{A1} \end{tabular}
& \begin{tabular}{@{}l@{}} $s_1$, $s_2$ \\ \texttt{B1},\ \texttt{A1} \end{tabular}
& \raisebox{1.5ex}{--}
& \begin{tabular}{@{}l@{}} $s_2$ \\ \texttt{B1} \end{tabular} \\
\hline
$2$
& \raisebox{1.5ex}{--}
& \raisebox{1.5ex}{--}
& \raisebox{1.5ex}{--}
& \raisebox{1.5ex}{--}
& \begin{tabular}{@{}l@{}} $s_1$ \\ \texttt{A3} \end{tabular}
& \raisebox{1.5ex}{--}
& \raisebox{1.5ex}{--}
& \begin{tabular}{@{}l@{}} $s_1$, $s_3$ \\ \texttt{B3},\ \texttt{A3} \end{tabular} \\
\hline
\end{tabular}
\end{center}
\end{table}

% === CB (modern/base style; C-block notation) ===
\begin{table}[ht!]
\begin{center}
\caption{{\bf Reachability \texttt{cB}.} This is the generic case where the block $B(i)$ meets $C(i-1)$, which implies $\beta(i-1)<\alpha(i-1)$. For the relevant $c$-factors, the distances to the end of the $C$-block are $d(\texttt{c4}) = 2m+\beta'$ and $d(\texttt{c8}) = \beta',$ with corresponding lengths $|\texttt{c4}| = k-\beta'$ and $|\texttt{c8}| = m.$} %In these tables we use the convention $\alpha(i-1)=\alpha', \beta(i-1)=\beta' $ and $\gamma(i-1)=\gamma'$.}
\label{tab:CB}
\vspace{-3pt}
\setlength{\tabcolsep}{3pt}
\renewcommand{\arraystretch}{1.15}
\begin{tabular}{|c|c|c|c|c|c|c|c|c|}
\hline
\texttt{cB}
& $\underset{m-\alpha}{\texttt{A1} \text{``1''}}$
& $\underset{\alpha}{\texttt{A2} \text{``0''}}$
& $\underset{m-\alpha}{\texttt{A3} \text{``2''}}$
& $\underset{\gamma}{\texttt{B1} \text{``1''}}$
& $\underset{m-\gamma}{\texttt{B2} \text{``3''}}$
& $\underset{\gamma}{\texttt{B3} \text{``2''}}$
& $\underset{m-\alpha}{\texttt{B4} \text{``0''}}$
& $\underset{\gamma}{\texttt{B5} \text{``3''}}$ \\
\hline\hline
$0$
& \begin{tabular}{@{}l@{}} $s_1$, $s_2$ \\ \texttt{c8},\ \texttt{c4} \end{tabular}
& \raisebox{1.5ex}{--}
& \begin{tabular}{@{}l@{}} $s_2$, $s_3$ \\ \texttt{c8},\ \texttt{c4} \end{tabular}
& \begin{tabular}{@{}l@{}} $s_1$, $s_2$ \\ \texttt{A2},\ \texttt{c8} \end{tabular}
& \begin{tabular}{@{}l@{}} $s_3$ \\ \texttt{c8} \end{tabular}
& \begin{tabular}{@{}l@{}} $s_2$, $s_3$ \\ \texttt{A2},\ \texttt{c8} \end{tabular}
& \raisebox{1.5ex}{--}
& \begin{tabular}{@{}l@{}} $s_1$, $s_3$ \\ \texttt{B4},\ \texttt{A2} \end{tabular} \\
\hline
$1$
& \raisebox{1.5ex}{--}
& \raisebox{1.5ex}{--}
& \begin{tabular}{@{}l@{}} $s_1$ \\ \texttt{A1} \end{tabular}
& \raisebox{1.5ex}{--}
& \begin{tabular}{@{}l@{}} $s_2$ \\ \texttt{A1} \end{tabular}
& \begin{tabular}{@{}l@{}} $s_1$, $s_2$ \\ \texttt{B1},\ \texttt{A1} \end{tabular}
& \raisebox{1.5ex}{--}
& \begin{tabular}{@{}l@{}} $s_2$ \\ \texttt{B1} \end{tabular} \\
\hline
$2$
& \raisebox{1.5ex}{--}
& \raisebox{1.5ex}{--}
& \raisebox{1.5ex}{--}
& \raisebox{1.5ex}{--}
& \begin{tabular}{@{}l@{}} $s_1$ \\ \texttt{A3} \end{tabular}
& \raisebox{1.5ex}{--}
& \raisebox{1.5ex}{--}
& \begin{tabular}{@{}l@{}} $s_1$, $s_3$ \\ \texttt{B3},\ \texttt{A3} \end{tabular} \\
\hline
\end{tabular}
\end{center}
\end{table}

% === BC (modern/base style; correct C-factor lengths) ===
\begin{table}[ht!]
\begin{center}
\caption{{\bf Reachability \texttt{bC}.} This is the generic case, where the block $C(i)$ meets $B(i-1)$, which implies $\beta(i)<\alpha(i)$. For the factor \texttt{b4}, the distance to the current $B$-block is $d(\texttt{b4}) = \gamma'$, and its length is $|\texttt{b4}|=m-\alpha'$.}
\label{tab:BC}
\vspace{-3pt}
\setlength{\tabcolsep}{3pt}
\renewcommand{\arraystretch}{1.15}
\begin{tabular}{|c|c|c|c|c|c|c|c|c|c|c|c|c|}
\hline
\texttt{bC}
& $\underset{m-\alpha}{\texttt{A1} \text{``1''}}$
& $\underset{\alpha}{\texttt{A2} \text{``0''}}$
& $\underset{m-\alpha}{\texttt{A3} \text{``2''}}$
& $\underset{m}{\texttt{C1} \text{``1''}}$
& $\underset{m}{\texttt{C2} \text{``2''}}$
& $\underset{\beta}{\texttt{C3} \text{``1''}}$
& $\underset{k-\beta}{\texttt{C4} \text{``0''}}$
& $\underset{m-k}{\texttt{C5} \text{``3''}}$
& $\underset{\beta}{\texttt{C6} \text{``2''}}$
& $\underset{k-\beta}{\texttt{C7} \text{``3''}}$
& $\underset{m}{\texttt{C8} \text{``0''}}$
& $\underset{\beta}{\texttt{C9} \text{``3''}}$ \\
\hline\hline
$0$
& \begin{tabular}{@{}l@{}} $s_1$ \\ \texttt{b4} \end{tabular}
& \raisebox{1.5ex}{--}
& \begin{tabular}{@{}l@{}} $s_2$ \\ \texttt{b4} \end{tabular}
& \begin{tabular}{@{}l@{}} $s_1$, $s_2$ \\ \texttt{A2},\ \texttt{b4} \end{tabular}
& \begin{tabular}{@{}l@{}} $s_2$, $s_3$ \\ \texttt{A2},\ \texttt{b4} \end{tabular}
& \begin{tabular}{@{}l@{}} $s_2$ \\ \texttt{A2} \end{tabular}
& \raisebox{1.5ex}{--}
& \begin{tabular}{@{}l@{}} $s_3$ \\ \texttt{A2} \end{tabular}
& \begin{tabular}{@{}l@{}} $s_3$ \\ \texttt{A2} \end{tabular}
& \begin{tabular}{@{}l@{}} $s_1$ \\ \texttt{C4} \end{tabular}
& \raisebox{1.5ex}{--}
& \begin{tabular}{@{}l@{}} $s_1$ \\ \texttt{C8} \end{tabular} \\
\hline
$1$
& \raisebox{1.5ex}{--}
& \raisebox{1.5ex}{--}
& \begin{tabular}{@{}l@{}} $s_1$ \\ \texttt{A1} \end{tabular}
& \raisebox{1.5ex}{--}
& \begin{tabular}{@{}l@{}} $s_1$, $s_2$ \\ \texttt{C1},\ \texttt{A1} \end{tabular}
& \raisebox{1.5ex}{--}
& \raisebox{1.5ex}{--}
& \begin{tabular}{@{}l@{}} $s_2$ \\ \texttt{C1} \end{tabular}
& \begin{tabular}{@{}l@{}} $s_1$, $s_2$ \\ \texttt{C3},\ \texttt{C1} \end{tabular}
& \begin{tabular}{@{}l@{}} $s_2$ \\ \texttt{C1} \end{tabular}
& \raisebox{1.5ex}{--}
& \begin{tabular}{@{}l@{}} $s_2$ \\ \texttt{C3} \end{tabular} \\
\hline
$2$
& \raisebox{1.5ex}{--}
& \raisebox{1.5ex}{--}
& \raisebox{1.5ex}{--}
& \raisebox{1.5ex}{--}
& \raisebox{1.5ex}{--}
& \raisebox{1.5ex}{--}
& \raisebox{1.5ex}{--}
& \begin{tabular}{@{}l@{}} $s_1$ \\ \texttt{C2} \end{tabular}
& \raisebox{1.5ex}{--}
& \begin{tabular}{@{}l@{}} $s_3$ \\ \texttt{A3} \end{tabular}
& \raisebox{1.5ex}{--}
& \begin{tabular}{@{}l@{}} $s_3$ \\ \texttt{C2} \end{tabular} \\
\hline
\end{tabular}
\end{center}
\end{table}

% === CC (modern/base style; C-block notation + correct C-lengths) ===
\begin{table}[ht!]
\begin{center}
\caption{{\bf Reachability \texttt{cC}.} This is the generic case, where the block $C(i)$ meets $C(i-1)$, which implies $\beta(i)<\alpha(i)$ and $\beta(i-1)<\alpha(i-1)$. For the relevant $c$-factors, the distances to the end of the $C$-block are $d(\texttt{c4}) = 2m+\beta'$ and $d(\texttt{c8}) = \beta',$ with corresponding lengths $|\texttt{c4}| = k-\beta'$ and $|\texttt{c8}| = m.$}
\label{tab:CC}
\vspace{-3pt}
\setlength{\tabcolsep}{3pt}
\renewcommand{\arraystretch}{1.15}
\begin{tabular}{|c|c|c|c|c|c|c|c|c|c|c|c|c|}
\hline
\texttt{cC}
& $\underset{m-\alpha}{\texttt{A1} \text{``1''}}$
& $\underset{\alpha}{\texttt{A2} \text{``0''}}$
& $\underset{m-\alpha}{\texttt{A3} \text{``2''}}$
& $\underset{m}{\texttt{C1} \text{``1''}}$
& $\underset{m}{\texttt{C2} \text{``2''}}$
& $\underset{\beta}{\texttt{C3} \text{``1''}}$
& $\underset{k-\beta}{\texttt{C4} \text{``0''}}$
& $\underset{m-k}{\texttt{C5} \text{``3''}}$
& $\underset{\beta}{\texttt{C6} \text{``2''}}$
& $\underset{k-\beta}{\texttt{C7} \text{``3''}}$
& $\underset{m}{\texttt{C8} \text{``0''}}$
& $\underset{\beta}{\texttt{C9} \text{``3''}}$ \\
\hline\hline
$0$
& \begin{tabular}{@{}l@{}} $s_1$, $s_2$ \\ \texttt{c8},\ \texttt{c4} \end{tabular}
& \raisebox{1.5ex}{--}
& \begin{tabular}{@{}l@{}} $s_2$, $s_3$ \\ \texttt{c8},\ \texttt{c4} \end{tabular}
& \begin{tabular}{@{}l@{}} $s_1$, $s_2$ \\ \texttt{A2},\ \texttt{c8} \end{tabular}
& \begin{tabular}{@{}l@{}} $s_2$, $s_3$ \\ \texttt{A2},\ \texttt{c8} \end{tabular}
& \begin{tabular}{@{}l@{}} $s_2$ \\ \texttt{A2} \end{tabular}
& \raisebox{1.5ex}{--}
& \begin{tabular}{@{}l@{}} $s_3$ \\ \texttt{A2} \end{tabular}
& \begin{tabular}{@{}l@{}} $s_3$ \\ \texttt{A2} \end{tabular}
& \begin{tabular}{@{}l@{}} $s_1$ \\ \texttt{C4} \end{tabular}
& \raisebox{1.5ex}{--}
& \begin{tabular}{@{}l@{}} $s_1$ \\ \texttt{C8} \end{tabular} \\
\hline
$1$
& \raisebox{1.5ex}{--}
& \raisebox{1.5ex}{--}
& \begin{tabular}{@{}l@{}} $s_1$ \\ \texttt{A1} \end{tabular}
& \raisebox{1.5ex}{--}
& \begin{tabular}{@{}l@{}} $s_1$, $s_2$ \\ \texttt{C1},\ \texttt{A1} \end{tabular}
& \raisebox{1.5ex}{--}
& \raisebox{1.5ex}{--}
& \begin{tabular}{@{}l@{}} $s_2$ \\ \texttt{C1} \end{tabular}
& \begin{tabular}{@{}l@{}} $s_1$, $s_2$ \\ \texttt{C3},\ \texttt{C1} \end{tabular}
& \begin{tabular}{@{}l@{}} $s_2$ \\ \texttt{C1} \end{tabular}
& \raisebox{1.5ex}{--}
& \begin{tabular}{@{}l@{}} $s_2$ \\ \texttt{C3} \end{tabular} \\
\hline
$2$
& \raisebox{1.5ex}{--}
& \raisebox{1.5ex}{--}
& \raisebox{1.5ex}{--}
& \raisebox{1.5ex}{--}
& \raisebox{1.5ex}{--}
& \raisebox{1.5ex}{--}
& \raisebox{1.5ex}{--}
& \begin{tabular}{@{}l@{}} $s_1$ \\ \texttt{C2} \end{tabular}
& \raisebox{1.5ex}{--}
& \begin{tabular}{@{}l@{}} $s_3$ \\ \texttt{A3} \end{tabular}
& \raisebox{1.5ex}{--}
& \begin{tabular}{@{}l@{}} $s_3$ \\ \texttt{C2} \end{tabular} \\
\hline
\end{tabular}
\end{center}
\end{table}

%%%%%%%%%%%Anti-collision%%%%%%%%%%%
\begin{table}[ht!]
\caption{{\bf Anti-collision \texttt{sink\!\! B}.} This table concerns the first $B$-block \texttt{sink\!\! B}. 
Each cell lists a triple of inequalities corresponding to the moves
$s_1$, $s_2$, and $s_3$, respectively. 
For example (with a fixed reference point),
\(\max(\texttt{A1})+s_1<\min(\texttt{B1})\),
whereas
\(\min(\texttt{A1})+s_>\max(\texttt{B1})\).
Subscripts indicate the lengths of the corresponding factors.
In each cell, the first occurrence of the inequality symbol ``$>$''
is tight.}\label{tab:ACsinkB}
\centering
%\begin{table}[ht!]
%\caption{{\bf Anticollision} table for system \texttt{sink\!\! B}. Cell entries are the triple of relations concerning the moves $s_1$, $s_2$ and $s_3$. For example (with fixed reference point) $\max\texttt{A1}+$s_1$<\min\texttt{B1}$, but  $\min\texttt{A1}+$s_2$>\max\texttt{B1}$. The subscripts take a note of the lengths of the factors. Note that, in each case, the first occurence of the inequalitity ``$>$'' is tight.}\label{tab:ACsinkB}
%\centering
\vspace{-3pt}
\setlength{\tabcolsep}{2pt}
\renewcommand{\arraystretch}{1.05}
%\noindent\textbf{Legend:} \(\begin{array}{l}$s_1$=m\\$s_2$=2m{+}k\\$s_3$=3m{+}k\end{array}\)
%\resizebox{\textwidth}{!}{%
\begin{tabular}{|c|c|c|c|c|c|c|c|c|}
\hline
%\multicolumn{9}{|c|}{Anticollision: \texttt{sink\!\! B}} \\
%\hline
\texttt{ sink\!\! B\; }& $\underset{m}{\texttt{A1} \text{``1''}}$
& $\underset{0}{\texttt{A2} \text{``0''}}$
& $\underset{m}{\texttt{A3} \text{``2''}}$
& $\underset{k}{\texttt{B1} \text{``1''}}$
& $\underset{m-k}{\texttt{B2} \text{``3''}}$
& $\underset{k}{\texttt{B3} \text{``2''}}$
& $\underset{m}{\texttt{B4} \text{``0''}}$
& $\underset{k}{\texttt{B5} \text{``3''}}$ \\
% & \texttt{A1} ``1'' & \texttt{A2} ``0'' & \texttt{A3} ``2'' & \texttt{B1} ``1'' & \texttt{B2} ``3'' & \texttt{B3} ``2'' & \texttt{B4} ``0'' & \texttt{B5} ``3''\\
\hline\hline
\multicolumn{9}{|l|}{\textrm{value 0}}\\
\hline
\raisebox{1.5ex}{\texttt{sink}} & \raisebox{1.5ex}{--} & \raisebox{1.5ex}{$\varnothing$} & \raisebox{1.5ex}{--} & \raisebox{1.5ex}{--} & \raisebox{1.5ex}{--} & \raisebox{1.5ex}{--} &  \shortstack{\rule{0pt}{1.6ex}$\small{<}$\\$\small{<}$\\$\small{<}$} & \raisebox{1.5ex}{--}\\
\hline
\multicolumn{9}{|l|}{\textrm{value 1}}\\
\hline
\raisebox{.6ex}{\shortstack{\texttt{A1}\\{\scriptsize $m$}}} & \raisebox{1.5ex}{--} & \raisebox{1.5ex}{--} & \raisebox{1.5ex}{--} & \shortstack{\rule{0pt}{1.6ex}$\small{<}$\\$\large{\boldsymbol{>}}$\\$\large{\boldsymbol{>}}$} & \raisebox{1.5ex}{--} & \raisebox{1.5ex}{--} & \raisebox{1.5ex}{--} & \raisebox{1.5ex}{--}\\
\hline
\multicolumn{9}{|l|}{\textrm{value 2}}\\
\hline
\raisebox{.6ex}{\shortstack{\texttt{A3}\\{\scriptsize $m$}}} & \raisebox{1.5ex}{--} & \raisebox{1.5ex}{--} & \raisebox{1.5ex}{--} & \raisebox{1.5ex}{--} & \raisebox{1.5ex}{--} & \shortstack{\rule{0pt}{1.6ex}$\small{<}$\\$\large{\boldsymbol{>}}$\\$\large{\boldsymbol{>}}$} & \raisebox{1.5ex}{--} & \raisebox{1.5ex}{--}\\
\hline
\multicolumn{9}{|l|}{\textrm{value 3}}\\
\hline
\raisebox{2.5ex}{$\underset{m-k}{\texttt{B2}}$}
%\shortstack{\texttt{B2}\\{\scriptsize{$m\!-\!k$}}} 
& \raisebox{1.5ex}{--} & \raisebox{1.5ex}{--} & \raisebox{1.5ex}{--} & \raisebox{1.5ex}{--} & \raisebox{1.5ex}{--} & \raisebox{1.5ex}{--} & \raisebox{1.5ex}{--} & \shortstack{\rule{0pt}{1.6ex}$\small{<}$\\$\large{\boldsymbol{>}}$\\$\large{\boldsymbol{>}}$}\\
\hline
\end{tabular}
\end{table}

\begin{table}[ht!]
\caption{{\bf Anti-collision \texttt{bZB}.} This table concerns the block concatenation
$B(i-1)\zeta(i-1)B(i)$, whenever $\alpha(i)=0$. Thus, \texttt{Z} denotes the set of positions induced by the factor $\zeta(i-1)=0^m$. 
Some row factors have distance greater than $3m+k$ 
from the factor under consideration in the column. In such cases, any
preceding row factor with the same game value is omitted. For example,
$\min(\texttt{Z}+3m+k)>\max(\texttt{B4})$, and hence the factor \texttt{b4}
of the same value may be omitted. Recall that factors belonging to the
current block are denoted by capital letters, while factors from the
preceding block are denoted by lowercase letters. 
Row headers list the factor name and, for factors in the preceding block,
the distance to the end of the block,  
together with the factor length; for factors in the current block, only
the factor length is displayed. Note that the factors \texttt{b2} and
\texttt{A2} are empty.}\label{tab:antiBZB}
\vspace{-3pt}
\centering
\setlength{\tabcolsep}{2pt}
\renewcommand{\arraystretch}{1.05}
%\noindent\textbf{Legend:} \(\begin{array}{l}$s_1$=m\\$s_2$=2m{+}k\\$s_3$=3m{+}k\end{array}\)
%\resizebox{\textwidth}{!}{%
\begin{tabular}{|c|c|c|c|c|c|c|c|c|c|}
\hline
%\multicolumn{10}{|c|}{Anticollision: \texttt{BZB}} \\
%\hline
\texttt{bZB}&$\underset{m}{\texttt{Z}\text{``0''}}$
& $\underset{m}{\texttt{A1} \text{``1''}}$
& $\underset{0}{\texttt{A2}}$
& $\underset{m}{\texttt{A3} \text{``2''}}$
& $\underset{k}{\texttt{B1} \text{``1''}}$
& $\underset{m-k}{\texttt{B2} \text{``3''}}$
& $\underset{k}{\texttt{B3} \text{``2''}}$
& $\underset{m}{\texttt{B4} \text{``0''}}$
& $\underset{k}{\texttt{B5} \text{``3''}}$ \\
%&$0^m$ & \texttt{A1} ``1'' & \texttt{A2} ``0'' & \texttt{A3} ``2'' & \texttt{B1} ``1'' & \texttt{B2} ``3'' & \texttt{B3} ``2'' & \texttt{B4} ``0'' & \texttt{B5} ``3''\\
\hline\hline
\multicolumn{10}{|l|}{\textrm{value 0}}\\
\hline
\raisebox{2.5ex}{$\underset{m,k}{\texttt{b4}}$} &\shortstack{\rule{0pt}{1.6ex}$\small{<}$\\$\large{\boldsymbol{>}}$\\$\large{\boldsymbol{>}}$}& \raisebox{1.5ex}{--} & \raisebox{1.5ex}{$\varnothing$} & \raisebox{1.5ex}{--} & \raisebox{1.5ex}{--} & \raisebox{1.5ex}{--} & \raisebox{1.5ex}{--} &  \shortstack{\rule{0pt}{1.6ex}$\small{<}$\\$\small{<}$\\$\small{<}$} & \raisebox{1.5ex}{--}\\
\hline
\raisebox{2.5ex}{$\underset{m}{\texttt{Z}}$} &\raisebox{1.5ex}{--}& \raisebox{1.5ex}{--} & \raisebox{1.5ex}{$\varnothing$} & \raisebox{1.5ex}{--} & \raisebox{1.5ex}{--} & \raisebox{1.5ex}{--} & \raisebox{1.5ex}{--} &  \shortstack{\rule{0pt}{1.6ex}$\small{<}$\\$\small{<}$\\$\small{<}$} & \raisebox{1.5ex}{--}\\
\hline

\multicolumn{10}{|l|}{\textrm{value 1}}\\
\hline
%\raisebox{2.5ex}{$\underset{5m+k,k}{\texttt{a1}}$} &&  \shortstack{\rule{0pt}{1.6ex}$\small{<}$\\$\small{<}$\\$\small{<}$} & \raisebox{1.5ex}{--} & \raisebox{1.5ex}{--} &  \shortstack{\rule{0pt}{1.6ex}$\small{<}$\\$\small{<}$\\$\small{<}$} & \raisebox{1.5ex}{--} & \raisebox{1.5ex}{--} & \raisebox{1.5ex}{--} & \raisebox{1.5ex}{--}\\
%\hline
\raisebox{2.5ex}{$\underset{2m+k,m}{\texttt{b1}}$} &\raisebox{1.5ex}{--}&  \shortstack{\rule{0pt}{1.6ex}$\small{<}$\\$\small{<}$\\$\small{<}$} & \raisebox{1.5ex}{--} & \raisebox{1.5ex}{--} &  \shortstack{\rule{0pt}{1.6ex}$\small{<}$\\$\small{<}$\\$\small{<}$} & \raisebox{1.5ex}{--} & \raisebox{1.5ex}{--} & \raisebox{1.5ex}{--} & \raisebox{1.5ex}{--}\\
\hline
\raisebox{2.5ex}{$\underset{m}{\texttt{A1}}$} &\raisebox{1.5ex}{--}& \raisebox{1.5ex}{--} & \raisebox{1.5ex}{--} & \raisebox{1.5ex}{--} & \shortstack{\rule{0pt}{1.6ex}$\small{<}$\\$\large{\boldsymbol{>}}$\\$\large{\boldsymbol{>}}$} & \raisebox{1.5ex}{--} & \raisebox{1.5ex}{--} & \raisebox{1.5ex}{--} & \raisebox{1.5ex}{--}\\
\hline

\multicolumn{10}{|l|}{\textrm{value 2}}\\
\hline
%\raisebox{2.5ex}{$\underset{4m+2k}{\texttt{a3}}$} && \raisebox{1.5ex}{--} & \raisebox{1.5ex}{--} &  \shortstack{\rule{0pt}{1.6ex}$\small{<}$\\$\small{<}$\\$\small{<}$} & \raisebox{1.5ex}{--} & \raisebox{1.5ex}{--} &  \shortstack{\rule{0pt}{1.6ex}$\small{<}$\\$\small{<}$\\$\small{<}$} & \raisebox{1.5ex}{--} & \raisebox{1.5ex}{--}\\
%\hline
\raisebox{2.5ex}{$\underset{m+k,m}{\texttt{b3}}$} &\raisebox{1.5ex}{--}& \raisebox{1.5ex}{--} & \raisebox{1.5ex}{--} &  \shortstack{\rule{0pt}{1.6ex}$\small{<}$\\$\small{<}$\\$\small{<}$} & \raisebox{1.5ex}{--} & \raisebox{1.5ex}{--} &  \shortstack{\rule{0pt}{1.6ex}$\small{<}$\\$\small{<}$\\$\small{<}$} & \raisebox{1.5ex}{--} & \raisebox{1.5ex}{--}\\
\hline
\raisebox{2.5ex}{$\underset{m}{\texttt{A3}}$} &\raisebox{1.5ex}{--}& \raisebox{1.5ex}{--} & \raisebox{1.5ex}{--} & \raisebox{1.5ex}{--} & \raisebox{1.5ex}{--} & \raisebox{1.5ex}{--} & \shortstack{\rule{0pt}{1.6ex}$\small{<}$\\$\large{\boldsymbol{>}}$\\$\large{\boldsymbol{>}}$} & \raisebox{1.5ex}{--} & \raisebox{1.5ex}{--}\\
\hline
\multicolumn{10}{|l|}{\textrm{value 3}}\\
\hline
\raisebox{2.5ex}{$\underset{0,m}{\texttt{b5}}$} &\raisebox{1.5ex}{--}& \raisebox{1.5ex}{--} & \raisebox{1.5ex}{--} & \raisebox{1.5ex}{--} & \raisebox{1.5ex}{--} &  \shortstack{\rule{0pt}{1.6ex}$\small{<}$\\$\small{<}$\\$\small{<}$} & \raisebox{1.5ex}{--} & \raisebox{1.5ex}{--} &  \shortstack{\rule{0pt}{1.6ex}$\small{<}$\\$\small{<}$\\$\small{<}$}\\
\hline
\raisebox{2.5ex}{$\underset{k}{\texttt{B2}}$} &\raisebox{1.5ex}{--}& \raisebox{1.5ex}{--} & \raisebox{1.5ex}{--} & \raisebox{1.5ex}{--} & \raisebox{1.5ex}{--} & \raisebox{1.5ex}{--} & \raisebox{1.5ex}{--} & \raisebox{1.5ex}{--} & \shortstack{\rule{0pt}{1.6ex}$\small{<}$\\$\large{\boldsymbol{>}}$\\$\large{\boldsymbol{>}}$}\\
\hline
\end{tabular}%}
\end{table}

%\clearpage

\begin{table}[ht!]\caption{{\bf Anti-collision \texttt{bB}.} This table concerns the general block concatenation \texttt{bB}. The factor sets \texttt{a1}, \texttt{a2}, \texttt{a3} and \texttt{b2} are omitted since they are of at least a block length to the first possible collision, respectively.}
\label{tab:antiBB}
\centering
\vspace{-3pt}
\begin{tabular}{|c|c|c|c|c|c|c|c|c|}
\hline
\texttt{bB}& $\underset{m-\alpha}{\texttt{A1} \text{``1''}}$
& $\underset{\alpha}{\texttt{A2} \text{``0''}}$
& $\underset{m-\alpha}{\texttt{A3} \text{``2''}}$
& $\underset{\gamma}{\texttt{B1} \text{``1''}}$
& $\underset{m-\gamma}{\texttt{B2} \text{``3''}}$
& $\underset{\gamma}{\texttt{B3} \text{``2''}}$
& $\underset{m-\alpha}{\texttt{B4} \text{``0''}}$
& $\underset{\gamma}{\texttt{B5} \text{``3''}}$ \\
\hline\hline
\multicolumn{9}{|l|}{\textrm{value 0}}\\
\hline
\raisebox{2.5ex}{$\underset{\gamma',m-\alpha'}{\texttt{b4}}$} & \raisebox{1.5ex}{--} & \shortstack{\rule{0pt}{1.6ex}$\small{<}$\\$\large{\boldsymbol{>}}$\\$\large{\boldsymbol{>}}$} & \raisebox{1.5ex}{--} & \raisebox{1.5ex}{--} & \raisebox{1.5ex}{--} & \raisebox{1.5ex}{--} & \shortstack{\rule{0pt}{1.6ex}$\small{<}$\\$\small{<}$\\$\small{<}$} & \raisebox{1.5ex}{--}\\
\hline
\raisebox{2.5ex}{$\underset{\alpha}{\texttt{A2}}$} & \raisebox{1.5ex}{--} & \raisebox{1.5ex}{--} & \raisebox{1.5ex}{--} & \raisebox{1.5ex}{--} & \raisebox{1.5ex}{--} & \raisebox{1.5ex}{--} & \shortstack{\rule{0pt}{1.6ex}$\small{<}$\\$\small{<}$\\$\large{\boldsymbol{>}}$} & \raisebox{1.5ex}{--}\\
\hline
\multicolumn{9}{|l|}{\textrm{value 1}}\\
%\hline
%\raisebox{2.5ex}{$\underset{3m+k+\gamma',m-\alpha'}{\texttt{a1}}$} & \shortstack{\rule{0pt}{1.6ex}$\small{<}$\\$\small{<}$\\$\small{<}$} & \raisebox{1.5ex}{--} & \raisebox{1.5ex}{--} & \shortstack{\rule{0pt}{1.6ex}$\small{<}$\\$\small{<}$\\$\small{<}$} & \raisebox{1.5ex}{--} & \raisebox{1.5ex}{--} & \raisebox{1.5ex}{--} & \raisebox{1.5ex}{--}\\
\hline
\raisebox{2.5ex}{$\underset{2m+k,\gamma'}{\texttt{b1}}$} & \shortstack{\rule{0pt}{1.6ex}$\small{<}$\\$\small{<}$\\$\large{\boldsymbol{>}}$} & \raisebox{1.5ex}{--} & \raisebox{1.5ex}{--} & \shortstack{\rule{0pt}{1.6ex}$\small{<}$\\$\small{<}$\\$\small{<}$} & \raisebox{1.5ex}{--} & \raisebox{1.5ex}{--} & \raisebox{1.5ex}{--} & \raisebox{1.5ex}{--}\\
\hline
\raisebox{2.5ex}{$\underset{m-\alpha}{\texttt{A1}}$} & \raisebox{1.5ex}{--} & \raisebox{1.5ex}{--} & \raisebox{1.5ex}{--} & \shortstack{\rule{0pt}{1.6ex}$\small{<}$\\$\large{\boldsymbol{>}}$\\$\large{\boldsymbol{>}}$} & \raisebox{1.5ex}{--} & \raisebox{1.5ex}{--} & \raisebox{1.5ex}{--} & \raisebox{1.5ex}{--}\\
\hline
\multicolumn{9}{|l|}{\textrm{value 2}}\\
\hline
%\raisebox{2.5ex}{$\underset{2m+k+\gamma',m-\alpha'}{\texttt{a3}}$} & \raisebox{1.5ex}{--} & \raisebox{1.5ex}{--} & \shortstack{\rule{0pt}{1.6ex}$\small{<}$\\$\small{<}$\\$\small{<}$} & \raisebox{1.5ex}{--} & \raisebox{1.5ex}{--} & \shortstack{\rule{0pt}{1.6ex}$\small{<}$\\$\small{<}$\\$\small{<}$} & \raisebox{1.5ex}{--} & \raisebox{1.5ex}{--}\\
%\hline
\raisebox{2.5ex}{$\underset{m+k,\gamma'}{\texttt{b3}}$} & \raisebox{1.5ex}{--} & \raisebox{1.5ex}{--} & \shortstack{\rule{0pt}{1.6ex}$\small{<}$\\$\small{<}$\\$\large{\boldsymbol{>}}$} & \raisebox{1.5ex}{--} & \raisebox{1.5ex}{--} & \shortstack{\rule{0pt}{1.6ex}$\small{<}$\\$\small{<}$\\$\small{<}$} & \raisebox{1.5ex}{--} & \raisebox{1.5ex}{--}\\
\hline
\raisebox{2.5ex}{$\underset{m-\alpha}{\texttt{A3}}$} & \raisebox{1.5ex}{--} & \raisebox{1.5ex}{--} & \raisebox{1.5ex}{--} & \raisebox{1.5ex}{--} & \raisebox{1.5ex}{--} & \shortstack{\rule{0pt}{1.6ex}$\small{<}$\\$\large{\boldsymbol{>}}$\\$\large{\boldsymbol{>}}$} & \raisebox{1.5ex}{--} & \raisebox{1.5ex}{--}\\
\hline
\multicolumn{9}{|l|}{\textrm{value 3}}\\
%\hline
%\raisebox{2.5ex}{$\underset{m+k+\gamma', m-\gamma'}{\texttt{b2}}$} & \raisebox{1.5ex}{--} & \raisebox{1.5ex}{--} & \raisebox{1.5ex}{--} & \raisebox{1.5ex}{--} & \shortstack{\rule{0pt}{1.6ex}$\small{<}$\\$\small{<}$\\$\small{<}$} & \raisebox{1.5ex}{--} & \raisebox{1.5ex}{--} & \shortstack{\rule{0pt}{1.6ex}$\small{<}$\\$\small{<}$\\$\small{<}$}\\
\hline
\raisebox{2.5ex}{$\underset{0,\gamma'}{\texttt{b5}}$} & \raisebox{1.5ex}{--} & \raisebox{1.5ex}{--} & \raisebox{1.5ex}{--} & \raisebox{1.5ex}{--} & \shortstack{\rule{0pt}{1.6ex}$\small{<}$\\$\small{<}$\\$\large{\boldsymbol{>}}$} & \raisebox{1.5ex}{--} & \raisebox{1.5ex}{--} & \shortstack{\rule{0pt}{1.6ex}$\small{<}$\\$\small{<}$\\$\small{<}$}\\
\hline 
\raisebox{2.5ex}{$\underset{m-\gamma}{\texttt{B2}}$} & \raisebox{1.5ex}{--} & \raisebox{1.5ex}{--} & \raisebox{1.5ex}{--} & \raisebox{1.5ex}{--} & \raisebox{1.5ex}{--} & \raisebox{1.5ex}{--} & \raisebox{1.5ex}{--} & \shortstack{\rule{0pt}{1.6ex}$\small{<}$\\$\large{\boldsymbol{>}}$\\$\large{\boldsymbol{>}}$}\\
\hline
\end{tabular}
\end{table}

\begin{table}[ht!]\caption{{\bf Anti-collision \texttt{cB}.} This table concerns the block concatenation \texttt{cB}. The factors of the $B$-block are omitted since they were already tested in Table~\ref{tab:antiBB}.}
\label{tab:antiCB}
\centering
\setlength{\tabcolsep}{2pt}
\renewcommand{\arraystretch}{1.05}
\vspace{-3pt}
\begin{tabular}{|c|c|c|c|c|c|c|c|c|}
\hline
\texttt{cB}& $\underset{m-\alpha}{\texttt{A1} \text{``1''}}$
& $\underset{\alpha}{\texttt{A2} \text{``0''}}$
& $\underset{m-\alpha}{\texttt{A3} \text{``2''}}$
& $\underset{\gamma}{\texttt{B1} \text{``1''}}$
& $\underset{m-\gamma}{\texttt{B2} \text{``3''}}$
& $\underset{\gamma}{\texttt{B3} \text{``2''}}$
& $\underset{m-\alpha}{\texttt{B4} \text{``0''}}$
& $\underset{\gamma}{\texttt{B5} \text{``3''}}$ \\
\hline\hline
\multicolumn{9}{|l|}{\textrm{value 0}}\\
%\hline
%\raisebox{2.5ex}{$\underset{5m+k+\beta',\alpha'}{\texttt{a2}}$} & \raisebox{1.5ex}{--} &  \shortstack{\rule{0pt}{1.6ex}$\small{<}$\\$\small{<}$\\$\small{<}$} & \raisebox{1.5ex}{--} & \raisebox{1.5ex}{--} & \raisebox{1.5ex}{--} & \raisebox{1.5ex}{--} &  \shortstack{\rule{0pt}{1.6ex}$\small{<}$\\$\small{<}$\\$\small{<}$} & \raisebox{1.5ex}{--}\\
\hline
\raisebox{2.5ex}{$\underset{2m+k,k-\beta'}{\texttt{c4}}$} & \raisebox{1.5ex}{--} & \shortstack{\rule{0pt}{1.6ex}$\small{<}$\\$\small{<}$\\$\large{\boldsymbol{>}}$} & \raisebox{1.5ex}{--} & \raisebox{1.5ex}{--} & \raisebox{1.5ex}{--} & \raisebox{1.5ex}{--} &  \shortstack{\rule{0pt}{1.6ex}$\small{<}$\\$\small{<}$\\$\small{<}$} & \raisebox{1.5ex}{--}\\
\hline
\raisebox{2.5ex}{$\underset{m+\beta',m}{\texttt{c8}}$} & \raisebox{1.5ex}{--} & \shortstack{\rule{0pt}{1.6ex}$\small{<}$\\$\large{\boldsymbol{>}}$\\$\large{\boldsymbol{>}}$} & \raisebox{1.5ex}{--} & \raisebox{1.5ex}{--} & \raisebox{1.5ex}{--} & \raisebox{1.5ex}{--} &  \shortstack{\rule{0pt}{1.6ex}$\small{<}$\\$\small{<}$\\$\small{<}$} & \raisebox{1.5ex}{--}\\
\hline
%\raisebox{2.5ex}{$\underset{\alpha}{\texttt{A2}}$} & \raisebox{1.5ex}{--} & \raisebox{1.5ex}{--} & \raisebox{1.5ex}{--} & \raisebox{1.5ex}{--} & \raisebox{1.5ex}{--} & \raisebox{1.5ex}{--} & \shortstack{\rule{0pt}{1.6ex}$\small{<}$\\$\small{<}$\\$\large{\boldsymbol{>}}$} & \raisebox{1.5ex}{--}\\
%\hline
\multicolumn{9}{|l|}{\textrm{value 1}}\\
%\hline
%\raisebox{2.5ex}{$\underset{4m+2k}{\texttt{a1}}$} &  \shortstack{\rule{0pt}{1.6ex}$\small{<}$\\$\small{<}$\\$\small{<}$} & \raisebox{1.5ex}{--} & \raisebox{1.5ex}{--} &  \shortstack{\rule{0pt}{1.6ex}$\small{<}$\\$\small{<}$\\$\small{<}$} & \raisebox{1.5ex}{--} & \raisebox{1.5ex}{--} & \raisebox{1.5ex}{--} & \raisebox{1.5ex}{--}\\
\hline
\raisebox{2.5ex}{$\underset{3m+k+\beta',m}{\texttt{c1}}$} &  \shortstack{\rule{0pt}{1.6ex}$\small{<}$\\$\small{<}$\\$\small{<}$} & \raisebox{1.5ex}{--} & \raisebox{1.5ex}{--} &  \shortstack{\rule{0pt}{1.6ex}$\small{<}$\\$\small{<}$\\$\small{<}$} & \raisebox{1.5ex}{--} & \raisebox{1.5ex}{--} & \raisebox{1.5ex}{--} & \raisebox{1.5ex}{--}\\
\hline
\raisebox{2.5ex}{$\underset{2m+k,\beta'}{\texttt{c3}}$} & \shortstack{\rule{0pt}{1.6ex}$\small{<}$\\$\small{<}$\\$\large{\boldsymbol{>}}$} & \raisebox{1.5ex}{--} & \raisebox{1.5ex}{--} &  \shortstack{\rule{0pt}{1.6ex}$\small{<}$\\$\small{<}$\\$\small{<}$} & \raisebox{1.5ex}{--} & \raisebox{1.5ex}{--} & \raisebox{1.5ex}{--} & \raisebox{1.5ex}{--}\\
\hline
%\raisebox{2.5ex}{$\underset{m-\alpha}{\texttt{A1}}$} & \raisebox{1.5ex}{--} & \raisebox{1.5ex}{--} & \raisebox{1.5ex}{--} & \shortstack{\rule{0pt}{1.6ex}$\small{<}$\\$\large{\boldsymbol{>}}$\\$\large{\boldsymbol{>}}$} & \raisebox{1.5ex}{--} & \raisebox{1.5ex}{--} & \raisebox{1.5ex}{--} & \raisebox{1.5ex}{--}\\
%\hline
\multicolumn{9}{|l|}{\textrm{value 2}}\\
%\hline
%\raisebox{2.5ex}{$\underset{5m+2k}{\texttt{a3}}$} & \raisebox{1.5ex}{--} & \raisebox{1.5ex}{--} &  \shortstack{\rule{0pt}{1.6ex}$\small{<}$\\$\small{<}$\\$\small{<}$} & \raisebox{1.5ex}{--} & \raisebox{1.5ex}{--} &  \shortstack{\rule{0pt}{1.6ex}$\small{<}$\\$\small{<}$\\$\small{<}$} & \raisebox{1.5ex}{--} & \raisebox{1.5ex}{--}\\
\hline
\raisebox{2.5ex}{$\underset{2m+k+\beta',m}{\texttt{c2}}$} & \raisebox{1.5ex}{--} & \raisebox{1.5ex}{--} &  \shortstack{\rule{0pt}{1.6ex}$\small{<}$\\$\small{<}$\\$\small{<}$} & \raisebox{1.5ex}{--} & \raisebox{1.5ex}{--} &  \shortstack{\rule{0pt}{1.6ex}$\small{<}$\\$\small{<}$\\$\small{<}$} & \raisebox{1.5ex}{--} & \raisebox{1.5ex}{--}\\
\hline
\raisebox{2.5ex}{$\underset{m+k,\beta'}{\texttt{c6}}$} & \raisebox{1.5ex}{--} & \raisebox{1.5ex}{--} & \shortstack{\rule{0pt}{1.6ex}$\small{<}$\\$\small{<}$\\$\large{\boldsymbol{>}}$} & \raisebox{1.5ex}{--} & \raisebox{1.5ex}{--} &  \shortstack{\rule{0pt}{1.6ex}$\small{<}$\\$\small{<}$\\$\small{<}$} & \raisebox{1.5ex}{--} & \raisebox{1.5ex}{--}\\
\hline
%\raisebox{2.5ex}{$\underset{m-\alpha}{\texttt{A3}}$} & \raisebox{1.5ex}{--} & \raisebox{1.5ex}{--} & \raisebox{1.5ex}{--} & \raisebox{1.5ex}{--} & \raisebox{1.5ex}{--} & \shortstack{\rule{0pt}{1.6ex}$\small{<}$\\$\large{\boldsymbol{>}}$\\$\large{\boldsymbol{>}}$} & \raisebox{1.5ex}{--} & \raisebox{1.5ex}{--}\\
%\hline
\multicolumn{9}{|l|}{\textrm{value 3}}\\
%\hline
%\raisebox{2.5ex}{$\underset{2m+\beta}{\texttt{c5}}$} & \raisebox{1.5ex}{--} & \raisebox{1.5ex}{--} & \raisebox{1.5ex}{--} & \raisebox{1.5ex}{--} & \shortstack{--\\--\\--} & \raisebox{1.5ex}{--} & \raisebox{1.5ex}{--} &  \shortstack{\rule{0pt}{1.6ex}$\small{<}$\\$\small{<}$\\$\small{<}$}\\
%\hline
%\raisebox{2.5ex}{$\underset{m+k}{\texttt{c7}}$} & \raisebox{1.5ex}{--} & \raisebox{1.5ex}{--} & \raisebox{1.5ex}{--} & \raisebox{1.5ex}{--} & \shortstack{--\\--\\--} & \raisebox{1.5ex}{--} & \raisebox{1.5ex}{--} &  \shortstack{\rule{0pt}{1.6ex}$\small{<}$\\$\small{<}$\\$\small{<}$}\\
\hline
\raisebox{2.5ex}{$\underset{0,\beta'}{\texttt{c9}}$} & \raisebox{1.5ex}{--} & \raisebox{1.5ex}{--} & \raisebox{1.5ex}{--} & \raisebox{1.5ex}{--} & \raisebox{1.5ex}{--} & \raisebox{1.5ex}{--} & \raisebox{1.5ex}{--} &  \shortstack{\rule{0pt}{1.6ex}$\small{<}$\\$\small{<}$\\$\small{<}$}\\
\hline
\end{tabular}
\end{table}

\begin{table}[ht!]\caption{{\bf Anti-collision \texttt{bC}.} This table concerns the block concatenation \texttt{bC}. The factor sets \texttt{a1}, \texttt{a2} and \texttt{a3} have been omitted since they are a block length apart their corresponding  values  in the $C$-block.%with $(m,k)=(7,6)$. Left block: $B(0)$ with $\alpha=0,\gamma=6$. Right block: $I(1)$ with $\alpha=6,\beta=5$. Cell entries are the triple of relations for ($s_1$,$s_2$,$s_3$) = $s+\langle G\rangle$ to $\langle F\rangle$.
}\label{tab:antiBC}
\centering
\setlength{\tabcolsep}{2pt}
\renewcommand{\arraystretch}{1.05}
%\noindent\textbf{Legend:} \(\begin{array}{l}$s_1$=m\\$s_2$=2m{+}k\\$s_3$=3m{+}k\end{array}\)
%\resizebox{\textwidth}{!}{%
\begin{tabular}{|c|c|c|c|c|c|c|c|c|c|c|c|c|}
\hline
%\multicolumn{13}{|c|}{\textbf{System \texttt{BC}}} \\
%\hline
\texttt{bC}
& $\underset{m-\alpha}{\texttt{A1} \text{``1''}}$
& $\underset{\alpha}{\texttt{A2} \text{``0''}}$
& $\underset{m-\alpha}{\texttt{A3} \text{``2''}}$
& $\underset{m}{\texttt{C1} \text{``1''}}$
& $\underset{m}{\texttt{C2} \text{``2''}}$
& $\underset{\beta}{\texttt{C3} \text{``1''}}$
& $\underset{k-\beta}{\texttt{C4} \text{``0''}}$
& $\underset{m-k}{\texttt{C5} \text{``3''}}$
& $\underset{\beta}{\texttt{C6} \text{``2''}}$
& $\underset{k-\beta}{\texttt{C7} \text{``3''}}$
& $\underset{m}{\texttt{C8} \text{``0''}}$
& $\underset{\beta}{\texttt{C9} \text{``3''}}$ \\
\hline\hline
\multicolumn{13}{|l|}{\textrm{value 0}}\\
\hline
\raisebox{2.5ex}{$\underset{\gamma', m-\alpha'}{\texttt{b4}}$} & \raisebox{1.5ex}{--} & \shortstack{\rule{0pt}{1.6ex}$\small{<}$\\$\large{\boldsymbol{>}}$\\$\large{\boldsymbol{>}}$} & \raisebox{1.5ex}{--} & \raisebox{1.5ex}{--} & \raisebox{1.5ex}{--} & \raisebox{1.5ex}{--} &  \shortstack{\rule{0pt}{1.6ex}$\small{<}$\\$\small{<}$\\$\small{<}$} & \raisebox{1.5ex}{--} & \raisebox{1.5ex}{--} & \raisebox{1.5ex}{--} &  \shortstack{\rule{0pt}{1.6ex}$\small{<}$\\$\small{<}$\\$\small{<}$} & \raisebox{1.5ex}{--}\\
\hline
\raisebox{2.5ex}{$\underset{\alpha}{\texttt{A2}}$} & \raisebox{1.5ex}{--} & \raisebox{1.5ex}{--} & \raisebox{1.5ex}{--} & \raisebox{1.5ex}{--} & \raisebox{1.5ex}{--} & \raisebox{1.5ex}{--} & \shortstack{\rule{0pt}{1.6ex}$\small{<}$\\$\small{<}$\\$\large{\boldsymbol{>}}$} & \raisebox{1.5ex}{--} & \raisebox{1.5ex}{--} & \raisebox{1.5ex}{--} &  \shortstack{\rule{0pt}{1.6ex}$\small{<}$\\$\small{<}$\\$\small{<}$} & \raisebox{1.5ex}{--}\\
\hline
\raisebox{2.5ex}{$\underset{k-\beta}{\texttt{C4}}$} & \raisebox{1.5ex}{--} & \raisebox{1.5ex}{--} & \raisebox{1.5ex}{--} & \raisebox{1.5ex}{--} & \raisebox{1.5ex}{--} & \raisebox{1.5ex}{--} & \raisebox{1.5ex}{--} & \raisebox{1.5ex}{--} & \raisebox{1.5ex}{--} & \raisebox{1.5ex}{--} & \shortstack{\rule{0pt}{1.6ex}$\small{<}$\\$\large{\boldsymbol{>}}$\\$\large{\boldsymbol{>}}$} & \raisebox{1.5ex}{--}\\
\hline
\multicolumn{13}{|l|}{\textrm{value 1}}\\
\hline
%\raisebox{2.5ex}{$\underset{3m+k+\gamma', m-\alpha'}{\texttt{a1}}$} &  \shortstack{\rule{0pt}{1.6ex}$\small{<}$\\$\small{<}$\\$\small{<}$} & \raisebox{1.5ex}{--} & \raisebox{1.5ex}{--} &  \shortstack{\rule{0pt}{1.6ex}$\small{<}$\\$\small{<}$\\$\small{<}$} & \raisebox{1.5ex}{--} &  \shortstack{\rule{0pt}{1.6ex}$\small{<}$\\$\small{<}$\\$\small{<}$} & \raisebox{1.5ex}{--} & \raisebox{1.5ex}{--} & \raisebox{1.5ex}{--} & \raisebox{1.5ex}{--} & \raisebox{1.5ex}{--} & \raisebox{1.5ex}{--}\\
%\hline
\raisebox{2.5ex}{$\underset{2m+k,\gamma'}{\texttt{b1}}$} & \shortstack{\rule{0pt}{1.6ex}$\small{<}$\\$\small{<}$\\$\large{\boldsymbol{>}}$} & \raisebox{1.5ex}{--} & \raisebox{1.5ex}{--} &  \shortstack{\rule{0pt}{1.6ex}$\small{<}$\\$\small{<}$\\$\small{<}$} & \raisebox{1.5ex}{--} &  \shortstack{\rule{0pt}{1.6ex}$\small{<}$\\$\small{<}$\\$\small{<}$} & \raisebox{1.5ex}{--} & \raisebox{1.5ex}{--} & \raisebox{1.5ex}{--} & \raisebox{1.5ex}{--} & \raisebox{1.5ex}{--} & \raisebox{1.5ex}{--}\\
\hline
\raisebox{2.5ex}{$\underset{m-\alpha}{\texttt{A1}}$} & \raisebox{1.5ex}{--} & \raisebox{1.5ex}{--} & \raisebox{1.5ex}{--} & \shortstack{\rule{0pt}{1.6ex}$\small{<}$\\$\large{\boldsymbol{>}}$\\$\large{\boldsymbol{>}}$} & \raisebox{1.5ex}{--} & \shortstack{\rule{0pt}{1.6ex}$\small{<}$\\$\small{<}$\\$\large{\boldsymbol{>}}$} & \raisebox{1.5ex}{--} & \raisebox{1.5ex}{--} & \raisebox{1.5ex}{--} & \raisebox{1.5ex}{--} & \raisebox{1.5ex}{--} & \raisebox{1.5ex}{--}\\
\hline
\raisebox{2.5ex}{$\underset{m}{\texttt{C1}}$} & \raisebox{1.5ex}{--} & \raisebox{1.5ex}{--} & \raisebox{1.5ex}{--} & \raisebox{1.5ex}{--} & \raisebox{1.5ex}{--} & \shortstack{\rule{0pt}{1.6ex}$\small{<}$\\$\large{\boldsymbol{>}}$\\$\large{\boldsymbol{>}}$} & \raisebox{1.5ex}{--} & \raisebox{1.5ex}{--} & \raisebox{1.5ex}{--} & \raisebox{1.5ex}{--} & \raisebox{1.5ex}{--} & \raisebox{1.5ex}{--}\\
\hline
\multicolumn{13}{|l|}{\textrm{value 2}}\\
\hline
%\raisebox{2.5ex}{$\underset{4m+k+\gamma',m-\alpha'}{\texttt{a3}}$} & \raisebox{1.5ex}{--} & \raisebox{1.5ex}{--} &  \shortstack{\rule{0pt}{1.6ex}$\small{<}$\\$\small{<}$\\$\small{<}$} & \raisebox{1.5ex}{--} &  \shortstack{\rule{0pt}{1.6ex}$\small{<}$\\$\small{<}$\\$\small{<}$} & \raisebox{1.5ex}{--} & \raisebox{1.5ex}{--} & \raisebox{1.5ex}{--} &  \shortstack{\rule{0pt}{1.6ex}$\small{<}$\\$\small{<}$\\$\small{<}$} & \raisebox{1.5ex}{--} & \raisebox{1.5ex}{--} & \raisebox{1.5ex}{--}\\
%\hline
\raisebox{2.5ex}{$\underset{m+k,\gamma'}{\texttt{b3}}$} & \raisebox{1.5ex}{--} & \raisebox{1.5ex}{--} & \shortstack{\rule{0pt}{1.6ex}$\small{<}$\\$\small{<}$\\$\large{\boldsymbol{>}}$} & \raisebox{1.5ex}{--} &  \shortstack{\rule{0pt}{1.6ex}$\small{<}$\\$\small{<}$\\$\small{<}$} & \raisebox{1.5ex}{--} & \raisebox{1.5ex}{--} & \raisebox{1.5ex}{--} &  \shortstack{\rule{0pt}{1.6ex}$\small{<}$\\$\small{<}$\\$\small{<}$} & \raisebox{1.5ex}{--} & \raisebox{1.5ex}{--} & \raisebox{1.5ex}{--}\\
\hline
\raisebox{2.5ex}{$\underset{m-\alpha}{\texttt{A3}}$} & \raisebox{1.5ex}{--} & \raisebox{1.5ex}{--} & \raisebox{1.5ex}{--} & \raisebox{1.5ex}{--} & \shortstack{\rule{0pt}{1.6ex}$\small{<}$\\$\large{\boldsymbol{>}}$\\$\large{\boldsymbol{>}}$} & \raisebox{1.5ex}{--} & \raisebox{1.5ex}{--} & \raisebox{1.5ex}{--} & \shortstack{\rule{0pt}{1.6ex}$\small{<}$\\$\small{<}$\\$\large{\boldsymbol{>}}$} & \raisebox{1.5ex}{--} & \raisebox{1.5ex}{--} & \raisebox{1.5ex}{--}\\
\hline
\raisebox{2.5ex}{$\underset{m}{\texttt{C2}}$} & \raisebox{1.5ex}{--} & \raisebox{1.5ex}{--} & \raisebox{1.5ex}{--} & \raisebox{1.5ex}{--} & \raisebox{1.5ex}{--} & \raisebox{1.5ex}{--} & \raisebox{1.5ex}{--} & \raisebox{1.5ex}{--} & \shortstack{\rule{0pt}{1.6ex}$\small{<}$\\$\large{\boldsymbol{>}}$\\$\large{\boldsymbol{>}}$} & \raisebox{1.5ex}{--} & \raisebox{1.5ex}{--} & \raisebox{1.5ex}{--}\\
\hline
\multicolumn{13}{|l|}{\textrm{value 3}}\\
%\hline
%\raisebox{2.5ex}{$\underset{2m+k}{\texttt{b2}}$} & \raisebox{1.5ex}{--} & \raisebox{1.5ex}{--} & \raisebox{1.5ex}{--} & \raisebox{1.5ex}{--} & \raisebox{1.5ex}{--} & \raisebox{1.5ex}{--} & \raisebox{1.5ex}{--} &  \shortstack{\rule{0pt}{1.6ex}$\small{<}$\\$\small{<}$\\$\small{<}$} & \raisebox{1.5ex}{--} &  \shortstack{\rule{0pt}{1.6ex}$\small{<}$\\$\small{<}$\\$\small{<}$} & \raisebox{1.5ex}{--} &  \shortstack{\rule{0pt}{1.6ex}$\small{<}$\\$\small{<}$\\$\small{<}$}\\
\hline
\raisebox{2.5ex}{$\underset{0,\gamma'}{\texttt{b5}}$} & \raisebox{1.5ex}{--} & \raisebox{1.5ex}{--} & \raisebox{1.5ex}{--} & \raisebox{1.5ex}{--} & \raisebox{1.5ex}{--} & \raisebox{1.5ex}{--} & \raisebox{1.5ex}{--} &  \shortstack{\rule{0pt}{1.6ex}$\small{<}$\\$\small{<}$\\$\small{<}$} & \raisebox{1.5ex}{--} &  \shortstack{\rule{0pt}{1.6ex}$\small{<}$\\$\small{<}$\\$\small{<}$} & \raisebox{1.5ex}{--} &  \shortstack{\rule{0pt}{1.6ex}$\small{<}$\\$\small{<}$\\$\small{<}$}\\
\hline
\raisebox{2.5ex}{$\underset{m-k}{\texttt{C5}}$} & \raisebox{1.5ex}{--} & \raisebox{1.5ex}{--} & \raisebox{1.5ex}{--} & \raisebox{1.5ex}{--} & \raisebox{1.5ex}{--} & \raisebox{1.5ex}{--} & \raisebox{1.5ex}{--} & \raisebox{1.5ex}{--} & \raisebox{1.5ex}{--} & \shortstack{\rule{0pt}{1.6ex}$\large{\boldsymbol{>}}$\\$\large{\boldsymbol{>}}$\\$\large{\boldsymbol{>}}$} & \raisebox{1.5ex}{--} & \shortstack{\rule{0pt}{1.6ex}$\small{<}$\\$\large{\boldsymbol{>}}$\\$\large{\boldsymbol{>}}$}\\
\hline
\raisebox{2.5ex}{$\underset{k-\beta}{\texttt{C7}}$} & \raisebox{1.5ex}{--} & \raisebox{1.5ex}{--} & \raisebox{1.5ex}{--} & \raisebox{1.5ex}{--} & \raisebox{1.5ex}{--} & \raisebox{1.5ex}{--} & \raisebox{1.5ex}{--} & \raisebox{1.5ex}{--} & \raisebox{1.5ex}{--} & \raisebox{1.5ex}{--} & \raisebox{1.5ex}{--} & \shortstack{\rule{0pt}{1.6ex}$\small{<}$\\$\large{\boldsymbol{>}}$\\$\large{\boldsymbol{>}}$}\\
\hline
\end{tabular}
\end{table}

%\clearpage

\begin{table}[ht!]\caption{{\bf Anti-collision \texttt{cC}.} This table concerns the block concatenation \texttt{cC}. The factors from the current $C$-factor are omitted since they were already tested in Table~\ref{tab:antiBC}.}
\label{tab:antiCC}
\vspace{-3pt}
\centering
\setlength{\tabcolsep}{2pt}
\renewcommand{\arraystretch}{1.05}
\begin{tabular}{|c|c|c|c|c|c|c|c|c|c|c|c|c|}
\hline
 \texttt{cC}& $\underset{m-\alpha}{\texttt{A1} \text{``1''}}$
& $\underset{\alpha}{\texttt{A2} \text{``0''}}$
& $\underset{m-\alpha}{\texttt{A3} \text{``2''}}$
& $\underset{m}{\texttt{C1} \text{``1''}}$
& $\underset{m}{\texttt{C2} \text{``2''}}$
& $\underset{\beta}{\texttt{C3} \text{``1''}}$
& $\underset{k-\beta}{\texttt{C4} \text{``0''}}$
& $\underset{m-k}{\texttt{C5} \text{``3''}}$
& $\underset{\beta}{\texttt{C6} \text{``2''}}$
& $\underset{k-\beta}{\texttt{C7} \text{``3''}}$
& $\underset{m}{\texttt{C8} \text{``0''}}$
& $\underset{\beta}{\texttt{C9} \text{``3''}}$ \\
\hline\hline
\multicolumn{13}{|l|}{\textrm{value 0}}\\
%\hline
%\raisebox{2.5ex}{$\underset{5m+k+\beta',\alpha'}{\texttt{a2}}$} & \raisebox{1.5ex}{--} &  \shortstack{\rule{0pt}{1.6ex}$\small{<}$\\$\small{<}$\\$\small{<}$} & \raisebox{1.5ex}{--} & \raisebox{1.5ex}{--} & \raisebox{1.5ex}{--} & \raisebox{1.5ex}{--} &  \shortstack{\rule{0pt}{1.6ex}$\small{<}$\\$\small{<}$\\$\small{<}$} & \raisebox{1.5ex}{--} & \raisebox{1.5ex}{--} & \raisebox{1.5ex}{--} &  \shortstack{\rule{0pt}{1.6ex}$\small{<}$\\$\small{<}$\\$\small{<}$} & \raisebox{1.5ex}{--}\\
\hline
\raisebox{2.5ex}{$\underset{2m+k,k-\beta'}{\texttt{c4}}$} & \raisebox{1.5ex}{--} & \shortstack{\rule{0pt}{1.6ex}$\small{<}$\\$\small{<}$\\$\large{\boldsymbol{>}}$} & \raisebox{1.5ex}{--} & \raisebox{1.5ex}{--} & \raisebox{1.5ex}{--} & \raisebox{1.5ex}{--} &  \shortstack{\rule{0pt}{1.6ex}$\small{<}$\\$\small{<}$\\$\small{<}$} & \raisebox{1.5ex}{--} & \raisebox{1.5ex}{--} & \raisebox{1.5ex}{--} &  \shortstack{\rule{0pt}{1.6ex}$\small{<}$\\$\small{<}$\\$\small{<}$} & \raisebox{1.5ex}{--}\\
\hline
\raisebox{2.5ex}{$\underset{m+\beta',m}{\texttt{c8}}$} & \raisebox{1.5ex}{--} & \shortstack{\rule{0pt}{1.6ex}$\small{<}$\\$\large{\boldsymbol{>}}$\\$\large{\boldsymbol{>}}$} & \raisebox{1.5ex}{--} & \raisebox{1.5ex}{--} & \raisebox{1.5ex}{--} & \raisebox{1.5ex}{--} &  \shortstack{\rule{0pt}{1.6ex}$\small{<}$\\$\small{<}$\\$\small{<}$} & \raisebox{1.5ex}{--} & \raisebox{1.5ex}{--} & \raisebox{1.5ex}{--} &  \shortstack{\rule{0pt}{1.6ex}$\small{<}$\\$\small{<}$\\$\small{<}$} & \raisebox{1.5ex}{--}\\
\hline
%\raisebox{2.5ex}{$\underset{\alpha}{\texttt{A2}}$} & \raisebox{1.5ex}{--} & \raisebox{1.5ex}{--} & \raisebox{1.5ex}{--} & \raisebox{1.5ex}{--} & \raisebox{1.5ex}{--} & \raisebox{1.5ex}{--} & \shortstack{\rule{0pt}{1.6ex}$\small{<}$\\$\small{<}$\\$\large{\boldsymbol{>}}$} & \raisebox{1.5ex}{--} & \raisebox{1.5ex}{--} & \raisebox{1.5ex}{--} &  \shortstack{\rule{0pt}{1.6ex}$\small{<}$\\$\small{<}$\\$\small{<}$} & \raisebox{1.5ex}{--}\\
%\hline
%\raisebox{2.5ex}{$\underset{k-\beta}{\texttt{C4}}$} & \raisebox{1.5ex}{--} & \raisebox{1.5ex}{--} & \raisebox{1.5ex}{--} & \raisebox{1.5ex}{--} & \raisebox{1.5ex}{--} & \raisebox{1.5ex}{--} & \raisebox{1.5ex}{--} & \raisebox{1.5ex}{--} & \raisebox{1.5ex}{--} & \raisebox{1.5ex}{--} & \shortstack{\rule{0pt}{1.6ex}$\small{<}$\\$\large{\boldsymbol{>}}$\\$\large{\boldsymbol{>}}$} & \raisebox{1.5ex}{--}\\
%\hline
\multicolumn{13}{|l|}{\textrm{value 1}}\\
\hline
\raisebox{2.5ex}{$\underset{3m+k+\beta',m}{\texttt{c1}}$} &  \shortstack{\rule{0pt}{1.6ex}$\small{<}$\\$\small{<}$\\$\small{<}$} & \raisebox{1.5ex}{--} & \raisebox{1.5ex}{--} &  \shortstack{\rule{0pt}{1.6ex}$\small{<}$\\$\small{<}$\\$\small{<}$} & \raisebox{1.5ex}{--} &  \shortstack{\rule{0pt}{1.6ex}$\small{<}$\\$\small{<}$\\$\small{<}$} & \raisebox{1.5ex}{--} & \raisebox{1.5ex}{--} & \raisebox{1.5ex}{--} & \raisebox{1.5ex}{--} & \raisebox{1.5ex}{--} & \raisebox{1.5ex}{--}\\
\hline
\raisebox{2.5ex}{$\underset{2m+k,\beta'}{\texttt{c3}}$} & \shortstack{\rule{0pt}{1.6ex}$\small{<}$\\$\small{<}$\\$\large{\boldsymbol{>}}$} & \raisebox{1.5ex}{--} & \raisebox{1.5ex}{--} &  \shortstack{\rule{0pt}{1.6ex}$\small{<}$\\$\small{<}$\\$\small{<}$} & \raisebox{1.5ex}{--} &  \shortstack{\rule{0pt}{1.6ex}$\small{<}$\\$\small{<}$\\$\small{<}$} & \raisebox{1.5ex}{--} & \raisebox{1.5ex}{--} & \raisebox{1.5ex}{--} & \raisebox{1.5ex}{--} & \raisebox{1.5ex}{--} & \raisebox{1.5ex}{--}\\
\hline
%\raisebox{2.5ex}{$\underset{m-\alpha}{\texttt{A1}}$} & \raisebox{1.5ex}{--} & \raisebox{1.5ex}{--} & \raisebox{1.5ex}{--} & \shortstack{\rule{0pt}{1.6ex}$\small{<}$\\$\large{\boldsymbol{>}}$\\$\large{\boldsymbol{>}}$} & \raisebox{1.5ex}{--} & \shortstack{\rule{0pt}{1.6ex}$\small{<}$\\$\small{<}$\\$\large{\boldsymbol{>}}$} & \raisebox{1.5ex}{--} & \raisebox{1.5ex}{--} & \raisebox{1.5ex}{--} & \raisebox{1.5ex}{--} & \raisebox{1.5ex}{--} & \raisebox{1.5ex}{--}\\
%\hline
%\raisebox{2.5ex}{$\underset{m}{\texttt{C1}}$} & \raisebox{1.5ex}{--} & \raisebox{1.5ex}{--} & \raisebox{1.5ex}{--} & \raisebox{1.5ex}{--} & \raisebox{1.5ex}{--} & \shortstack{\rule{0pt}{1.6ex}$\small{<}$\\$\large{\boldsymbol{>}}$\\$\large{\boldsymbol{>}}$} & \raisebox{1.5ex}{--} & \raisebox{1.5ex}{--} & \raisebox{1.5ex}{--} & \raisebox{1.5ex}{--} & \raisebox{1.5ex}{--} & \raisebox{1.5ex}{--}\\
%\hline
\multicolumn{13}{|l|}{\textrm{value 2}}\\
\hline
\raisebox{2.5ex}{$\underset{2m+k+\beta',m}{\texttt{c2}}$} & \raisebox{1.5ex}{--} & \raisebox{1.5ex}{--} &  \shortstack{\rule{0pt}{1.6ex}$\small{<}$\\$\small{<}$\\$\small{<}$} & \raisebox{1.5ex}{--} &  \shortstack{\rule{0pt}{1.6ex}$\small{<}$\\$\small{<}$\\$\small{<}$} & \raisebox{1.5ex}{--} & \raisebox{1.5ex}{--} & \raisebox{1.5ex}{--} &  \shortstack{\rule{0pt}{1.6ex}$\small{<}$\\$\small{<}$\\$\small{<}$} & \raisebox{1.5ex}{--} & \raisebox{1.5ex}{--} & \raisebox{1.5ex}{--}\\
\hline
\raisebox{2.5ex}{$\underset{m+k,\beta'}{\texttt{c6}}$} & \raisebox{1.5ex}{--} & \raisebox{1.5ex}{--} & \shortstack{\rule{0pt}{1.6ex}$\small{<}$\\$\small{<}$\\$\large{\boldsymbol{>}}$} & \raisebox{1.5ex}{--} &  \shortstack{\rule{0pt}{1.6ex}$\small{<}$\\$\small{<}$\\$\small{<}$} & \raisebox{1.5ex}{--} & \raisebox{1.5ex}{--} & \raisebox{1.5ex}{--} &  \shortstack{\rule{0pt}{1.6ex}$\small{<}$\\$\small{<}$\\$\small{<}$} & \raisebox{1.5ex}{--} & \raisebox{1.5ex}{--} & \raisebox{1.5ex}{--}\\
\hline
%\raisebox{2.5ex}{$\underset{m-\alpha}{\texttt{A3}}$} & \raisebox{1.5ex}{--} & \raisebox{1.5ex}{--} & \raisebox{1.5ex}{--} & \raisebox{1.5ex}{--} & \shortstack{\rule{0pt}{1.6ex}$\small{<}$\\$\large{\boldsymbol{>}}$\\$\large{\boldsymbol{>}}$} & \raisebox{1.5ex}{--} & \raisebox{1.5ex}{--} & \raisebox{1.5ex}{--} & \shortstack{\rule{0pt}{1.6ex}$\small{<}$\\$\small{<}$\\$\large{\boldsymbol{>}}$} & \raisebox{1.5ex}{--} & \raisebox{1.5ex}{--} & \raisebox{1.5ex}{--}\\
%\hline
%\raisebox{2.5ex}{$\underset{m}{\texttt{C2}}$} & \raisebox{1.5ex}{--} & \raisebox{1.5ex}{--} & \raisebox{1.5ex}{--} & \raisebox{1.5ex}{--} & \raisebox{1.5ex}{--} & \raisebox{1.5ex}{--} & \raisebox{1.5ex}{--} & \raisebox{1.5ex}{--} & \shortstack{\rule{0pt}{1.6ex}$\small{<}$\\$\large{\boldsymbol{>}}$\\$\large{\boldsymbol{>}}$} & \raisebox{1.5ex}{--} & \raisebox{1.5ex}{--} & \raisebox{1.5ex}{--}\\
%\hline
\multicolumn{13}{|l|}{\textrm{value 3}}\\
\hline
\raisebox{2.5ex}{$\underset{0,\beta'}{\texttt{c9}}$} & \raisebox{1.5ex}{--} & \raisebox{1.5ex}{--} & \raisebox{1.5ex}{--} & \raisebox{1.5ex}{--} & \raisebox{1.5ex}{--} & \raisebox{1.5ex}{--} & \raisebox{1.5ex}{--} &  \shortstack{\rule{0pt}{1.6ex}$\small{<}$\\$\small{<}$\\$\small{<}$} & \raisebox{1.5ex}{--} &  \shortstack{\rule{0pt}{1.6ex}$\small{<}$\\$\small{<}$\\$\small{<}$} & \raisebox{1.5ex}{--} &  \shortstack{\rule{0pt}{1.6ex}$\small{<}$\\$\small{<}$\\$\small{<}$}\\
\hline
%\raisebox{2.5ex}{$\underset{m-k}{\texttt{C5}}$} & \raisebox{1.5ex}{--} & \raisebox{1.5ex}{--} & \raisebox{1.5ex}{--} & \raisebox{1.5ex}{--} & \raisebox{1.5ex}{--} & \raisebox{1.5ex}{--} & \raisebox{1.5ex}{--} & \raisebox{1.5ex}{--} & \raisebox{1.5ex}{--} & \shortstack{$\large{\boldsymbol{>}}$\\$\large{\boldsymbol{>}}$\\$\large{\boldsymbol{>}}$} & \raisebox{1.5ex}{--} & \shortstack{\rule{0pt}{1.6ex}$\small{<}$\\$\large{\boldsymbol{>}}$\\$\large{\boldsymbol{>}}$}\\
%\hline
%\raisebox{2.5ex}{$\underset{k-\beta}{\texttt{C7}}$} & \raisebox{1.5ex}{--} & \raisebox{1.5ex}{--} & \raisebox{1.5ex}{--} & \raisebox{1.5ex}{--} & \raisebox{1.5ex}{--} & \raisebox{1.5ex}{--} & \raisebox{1.5ex}{--} & \raisebox{1.5ex}{--} & \raisebox{1.5ex}{--} & \raisebox{1.5ex}{--} & \raisebox{1.5ex}{--} & \shortstack{\rule{0pt}{1.6ex}$\small{<}$\\$\large{\boldsymbol{>}}$\\$\large{\boldsymbol{>}}$}\\
%\hline
\end{tabular}
\end{table}

\begin{figure}
     \centering
 \includegraphics[width=.45\linewidth]{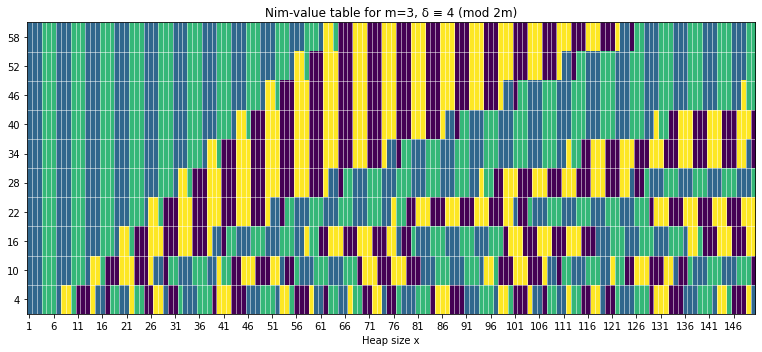}
 \includegraphics[width=.45\linewidth]{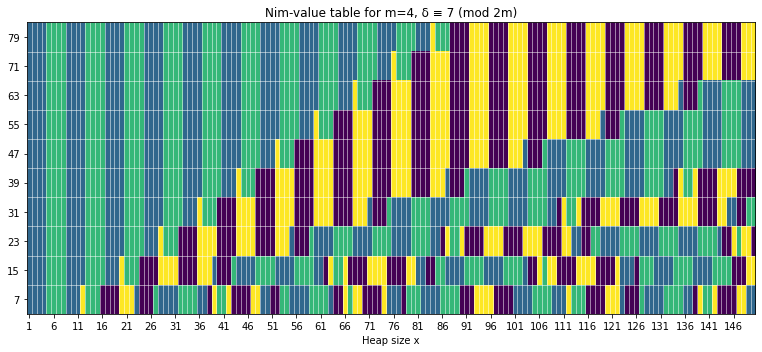}
     \caption{Nimbers for $m=3$, $\delta\equiv 4\pmod 6$ (left) and $m=4$ with  $\delta\equiv 7\pmod 8$ (right).}
     \label{fig:m3d4_m4d7}
 \end{figure}

\clearpage
It remains to justify the general case of item (ii), for a fixed $k$, and when $\delta = 2nm+k$, for $n\ge 1$. For each increase in $n$, the period length increases by $4m^2/{\gcd(m,\delta)}$, and this is due to the fact that each block is increased by exactly $4m$ positions via the two concatenated factors $1^m2^m$ and $3^m0^m$. Precisely, for each increase in $n$, insert $1^m2^m$ just after each $A$-prefix and insert $3^m0^m$ exactly between the last zero-factor and the last three-factor (independently of whether a $B$- or $C$-block); see Figure~\ref{fig:m3d4_m4d7}. The mex-justification is similar to case (i), but here it invokes another quadratic term. Apart from this observation, the mex verification remains the same as for the case $n=1$. This concludes the proof of Theorem~\ref{thm:main}.
\section{Final remarks}\label{sec:final}
On a more philosophical matter, with respect to the simplicity of the proof of Theorem~\ref{thm:folk}, one might ponder the value of understanding the fine details of a small sub-family of subtraction games. However this type of result has been on the list of conjectures for a long time, and our conjectured reciprocity with the wall convention requires further research. For example, Richard Guy writes in the first problem, in the first edition of {\em Unsolved problems in combinatorial games} \cite{guy1995unsolved} ``Subtraction games are known to be periodic. Investigate the relationship between the subtraction set and the length and structure of the period.''; see also \cite{nowakowski2019unsolved}. The second part of the conjecture for three move games \cite{Flammenkamp_1997, ward2016conjecture} concerns the non-additive three-move rulesets. Although they all are conjectured to have period length ``a sum of two moves'', and any such instance would individually be easily verified, the full classification appears to satisfy some fractal structure \cite{Flammenkamp_1997}, which has not yet been fully understood. See also the recent survey on subtraction games \cite{LarssonGoNC6}.

There is a closed formula expression for the $\mathcal P$-positions of {\sc additive wall subtraction} in \cite{berlekamp2004winning, Flammenkamp_1997}; and hence, by \eqref{eq:outwallsink}, it applies also to those of {\sc additive sink subtraction}. Namely the set of {\sc wall} $\mathcal P$-positions is generated by, for all $n\ge 0$,
$$n+\left\lfloor \frac{n}{m}\right\rfloor m + \left\lfloor\frac{2n}{2m+\delta - d} \right\rfloor(m+\delta)\ \text{ or }\ n+\left\lfloor \frac{n}{m}\right\rfloor m+ \left\lfloor\frac{2n}{\delta + d} \right\rfloor(m+\delta),$$ 
depending on whether $d\le m$ or not. To emphasize that the first (second) expression gives a linear (quadratic) period length, we rewrite this as
$$w_n=n+\left\lfloor \frac{n}{m}\right\rfloor m+ \left\lfloor\frac{2n}{km} \right\rfloor(m+\delta) \ \text{ or }\  w_n=n+\left\lfloor \frac{n}{m}\right\rfloor m+ \left\lfloor\frac{2n}{\ell m + 2d} \right\rfloor(m+\delta),$$ 
with $k = (2m+\delta-d)/m\in \mathbb N$ and $\ell = (\delta-d)/m\in \mathbb N$, respectively.  To illustrate, we get the following sequences for the parametrizations $m=2,\delta =d=1$ and $m=2,\delta =d=3$, respectively. 

\begin{center}
\begin{tabular}{c|cccccccccccc}
$n$   & 0 & 1 & 2 & 3 & 4  & 5  & 6  & 7  & 8  & 9  & 10 & 11 \\ \hline
$w_n(2,1)$ & 0 & 1 & 7 & 8 & 14 & 15 & 21 & 22 & 28 & 29 & 35 & 36\\ \hline
$s_n(2,1)$ & 5 & 6 & 12 & 13 & 19 & 20 & 26 & 27 & 33 & 34 & 40 & 41
\end{tabular}
\end{center}
\begin{center}
\begin{tabular}{c|cccccccccccc}
$n$   & 0 & 1 & 2 & 3  & 4  & 5  & 6  & 7  & 8  & 9  & 10 & 11 \\ \hline
$w_n(2,3)$ & 0 & 1 & 4 & 10 & 13 & 14 & 22 & 23 & 26 & 32 & 35 & 36 \\ \hline
$s_n(2,3)$ & 7 & 8 & 11 & 17 & 20 & 21 & 29 & 30 & 33 & 39 & 42 & 43
\end{tabular}
\end{center}
Here, we have included, as a second data row, the corresponding {\sc sink}-shifted $\mathcal P$-positions. 

Let us use the conjecture in Equation ~\ref{eq:conj} to directly compute all the nimbers of both {\sc sink} and {\sc wall} for the first ruleset, $S=\{2,3,5\}$ using also Ferguson's pairing principle. By the conjecture, from the $\mathcal P$-positions, we derive the nimber two positions for both conventions. 

For all $x\le 12$, if $v_w(x)=0$, then $v_s(14-3-x)=2$. Hence we get the {\sc sink} $2$-sequence: $3,4,10,11,\ldots$ 

Similarly, for all $x\le 17$, if $v_s(x)=0$, then $v_w(14+3-x)=2$. Hence we get the {\sc wall} $2$-sequence: $4,5,11,12,\ldots$ 

Next, the nimber one positions in the {\sc wall} convention are obtained via the classic Ferguson's pairing principle as 
$2,3,9,10,\ldots$. Now, the nimber three positions in the {\sc wall} convention must be the remaining positions, namely:  $6, 13, 20, \ldots$  

Similarly in the {\sc sink} convention, the conjecture gives the postions of the nimber ones by the formula $14-3-x$, where position $v_w(x)=1$. Thus, the sequence starts: $1,2,8,9, \ldots $ And the nimber three positions in the {\sc sink} convention must be the remaining positions, namely:  $5, 12, 19, \ldots$  

Note that all these values are obtained without any mex computation, from the classic $\mathcal P$-position formula in the {\sc wall} convention. Similarly one can fill in the periodic nimber positions for the ruleset $\{2,5,7\}$, without using any mex-rule. 

Let us summarize these findings. Let $S=\{s_1,s_2,s_3\}$ be an additive ruleset, with $s_1<s_2<s3$. 
\begin{enumerate}
\item By {\sc wall}-{\sc sink} outcome translation, $v_w(x)=0$ if and only if $v_s(x+s_3+1)=0$. 
\item By Ferguson, {\sc wall} zeros and ones are paired by the smallest move: $v_w(x)=0$ if and only $v_w(x+s_1)=1$. 
\item By our conjecture, given a {\sc wall} value  at position $x$, then $\sigma(v_w(x)) = v_s(-s_2-x)$. 
\item By our conjecture, given a {\sc sink} value  at position $x$, then $\sigma(v_s(x)) = v_w(s_2-x)$. 
\item Altogether by using the known formula of the $\mathcal P$-positions for the {\sc wall} convention, we can therefore compute the positions of all nimbers in $\{0,1,2\}$ for both conventions.
\item Since the ruleset has exactly three moves, we may conclude that, for each convention, the remaining positions have value three, respectively. 
\end{enumerate}

%Let us define the functions corresponding to all these nimbers.
From our duality conjecture together with Ferguson's pairing principle, let us derive explicit descriptions of all nimber classes in both conventions in terms of the wall $\mathcal P$-position sequence $(w_n)$ and the period $p$ of the outcome patterns.

For the {\sc wall} convention we get the sets of positions with nimbers smaller than three as:
\begin{align*}
\mathcal W_0 
  &= \{\, w_n  : n \ge 0 \,\}, \\[4pt]
\mathcal W_1 
  &= \{\, w_n + s_1 : n \ge 0 \,\}, \\[4pt]
\mathcal W_2 
  &= \{\, [\, - w_n - s_1 \,]_p + kp 
      : n \ge 0,\ k \ge 0 \,\}.
\end{align*}

And similarly for the {\sc sink} convention: 
\begin{align*}
\mathcal S_0 
  &= \{\, w_n + s_3+1 : n \ge 0 \,\}, \\[4pt]
\mathcal S_1 
  &= \{\, [\, -s_3 - w_n \,]_p + kp 
      : n \ge 0,\ k \ge 0 \,\}, \\[4pt]
\mathcal S_2 
  &= \{\, [\, -s_2 - w_n \,]_p + kp 
      : n \ge 0,\ k \ge 0 \,\}.
\end{align*}

As before, $[x]_p$ denotes the smallest nonnegative representative of $x$ modulo $p$. 
In either case the nimber three positions are the remaining ones.  Given our conjecture, in the case of ``primitive quadratic'' the nimber classes are also obtained from the {\sc wall} bracket formulas \cite{LM} by linear transformations and reduction modulo $p$.\\

%\[
%t_n 
% = n + 3m + \bigl(p-\delta\bigr)\left\lfloor\frac{n}{\delta}\right\rfloor,
% \qquad n\ge 0,
%\]

\noindent{\bf Acknowledgements.} We thank Richard Nowakowski for pointing out the reference \cite{althofer1995superlinear}. The images in Figure~\ref{fig:sinkwall} were generated with assistance from ChatGPT 5.2.
\bibliographystyle{plain}
\bibliography{references.bib}
\end{document}